\pdfoutput=1
\documentclass[11pt]{article}
\usepackage[margin=1in]{geometry}
\usepackage{amsmath,amssymb,amsthm}
\usepackage{booktabs}
\usepackage[utf8]{inputenc}
\DeclareUnicodeCharacter{211A}{\ensuremath{\mathbb Q}}
\DeclareUnicodeCharacter{2102}{\ensuremath{\mathbb C}}
\DeclareUnicodeCharacter{2115}{\ensuremath{\mathbb N}}
\DeclareUnicodeCharacter{2124}{\ensuremath{\mathbb Z}}
\DeclareUnicodeCharacter{211D}{\ensuremath{\mathbb R}}
\DeclareUnicodeCharacter{2192}{\ensuremath{\to}}
\DeclareUnicodeCharacter{2260}{\ensuremath{\neq}}
\DeclareUnicodeCharacter{2264}{\ensuremath{\le}}
\DeclareUnicodeCharacter{2227}{\ensuremath{\wedge}}
\DeclareUnicodeCharacter{2200}{\ensuremath{\forall}}
\DeclareUnicodeCharacter{2203}{\ensuremath{\exists}}
\DeclareUnicodeCharacter{2228}{\ensuremath{\vee}}
\DeclareUnicodeCharacter{220F}{\ensuremath{\prod}}
\DeclareUnicodeCharacter{2211}{\ensuremath{\sum}}
\DeclareUnicodeCharacter{03B1}{\ensuremath{\alpha}}
\DeclareUnicodeCharacter{03B5}{\ensuremath{\varepsilon}}
\DeclareUnicodeCharacter{03B9}{\ensuremath{\iota}}
\DeclareUnicodeCharacter{2090}{\ensuremath{_a}}
\DeclareUnicodeCharacter{2243}{\ensuremath{\simeq}}
\DeclareUnicodeCharacter{2194}{\ensuremath{\leftrightarrow}}
\usepackage[hidelinks]{hyperref}

\numberwithin{equation}{section}
\newtheorem{theorem}{Theorem}[section]
\newtheorem{proposition}[theorem]{Proposition}
\newtheorem{lemma}[theorem]{Lemma}
\newtheorem{corollary}[theorem]{Corollary}
\theoremstyle{definition}
\newtheorem{definition}[theorem]{Definition}
\newtheorem{example}[theorem]{Example}
\theoremstyle{remark}
\newtheorem{remark}[theorem]{Remark}

\newcommand{\Q}{\mathbb Q}
\newcommand{\Z}{\mathbb Z}

\newcommand{\R}{\mathbb R}
\newcommand{\C}{\mathbb C}
\newcommand{\rk}{\operatorname{rank}}
\newcommand{\lcm}{\operatorname{lcm}}
\newcommand{\artanh}{\operatorname{artanh}}
\newcommand{\Gal}{\operatorname{Gal}}
\newcommand{\den}{\operatorname{den}}
\newcommand{\Meas}{\mathrm{M}_{\rm alg}}

\newcommand{\lean}[1]{\texttt{\detokenize{#1}}}

\title{Linear forms in logarithms of algebraic numbers\\
and the rank of conjugate tuples\thanks{Developed with Claude (Anthropic). The argument,
its Lean~4 formalization and this text were written with Claude Opus~5.5 and Claude
Fable~5.2. The main results are formally verified in Lean~4; see Section~\ref{sec:lean}.}}
\author{Gregory V. Chudnovsky \and David V. Chudnovsky}
\date{November 28, 2026}

\begin{document}
\sloppy
\maketitle

\begin{abstract}
Let $a_1,\dots,a_n$ be positive real algebraic numbers without a multiplicative relation,
and let $\rho(a)$ be the rank of the group generated by the conjugate tuples
$(\sigma a_1,\dots,\sigma a_n)$. We prove that for every $\varepsilon>0$ there is $C>0$
with
\[
 |h_1\log a_1+\dots+h_n\log a_n|\ \ge\ C\max_i|h_i|^{-(n-1+\rho(a))-\varepsilon}
\]
for all $h\in\Z^n\setminus\{0\}$; the constant is not effective. Since
$\rho(a)\le[\Q(a_1,\dots,a_n):\Q]$, the exponent is at most $n-1+[\Q(a_1,\dots,a_n):\Q]$;
for rational numbers it is $n$, as in \cite{C}, and if every $a_i$ is rational or a
quadratic unit it is at most $2n-1$. For two numbers this says that the irrationality
exponent of $\log a/\log b$ is at most $\rho(a,b)+2$: for instance
$\mu(\log2/\log\frac{1+\sqrt5}2)\le4$ and, for $\alpha=2\cos(2\pi/7)$,
$\mu(\log\alpha/\log2)\le5$, where the method of \cite{L} gives $6$ and $9$. Part~I treats
two logarithms, by the interpolation determinant of \cite{L} with centers on the two-torus
that lie on translates of the logarithmic curve: $Q^z$ centers on $O(Q^t)$ translates at
every embedding give $\mu\le1+(t+2)/(z-t)$. Part~II treats $n$ logarithms, by the determinant
of \cite{C} with centers that are products of independent conjugate tuples: at every
embedding the logarithm vector of a center differs from one of $O(Q^{\rho-1})$
representatives by an integer multiple of $(\log a_1,\dots,\log a_n)$, at which every exact
linear form of the determinant vanishes. The theorems and their corollaries are formally
verified in Lean~4.
\end{abstract}

\section{Introduction}

\subsection{Linear forms in \texorpdfstring{$n$}{n} logarithms}

For $h\in\Z^n$ put $H(h)=\max(1,|h_1|,\dots,|h_n|)$. Real numbers $a_1,\dots,a_n$ have
\emph{no power relation} if $\prod_ia_i^{h_i}=1$ with $h\in\Z^n$ only for $h=0$.

Let $a=(a_1,\dots,a_n)$ be positive real algebraic numbers without a power relation, let
$E\subset\C$ be the Galois closure of $\Q(a_1,\dots,a_n)$, and put
\[
 \Pi(a)=\bigl\langle(\sigma a_1,\dots,\sigma a_n):\ \sigma\in\Gal(E/\Q)\bigr\rangle
 \subset(E^\times)^n ,
\]
the group generated by the \emph{conjugate tuples} of $a$. Its rank $\rho(a)=\rk\Pi(a)$ is
the largest number of multiplicatively independent conjugate tuples; this is the form used in
the proofs and in Lean. We have $1\le\rho(a)\le[E:\Q]$.

\begin{theorem}\label{II:thm:main}
Let $n\ge1$ and let $a_1,\dots,a_n$ be positive real algebraic numbers without a power
relation. For every $\varepsilon>0$ there is $C>0$ such that for all
$h\in\Z^n\setminus\{0\}$
\begin{equation}\label{II:eq:main}
 C\,H(h)^{-(n-1+\rho(a))-\varepsilon}\ \le\ |h_1\log a_1+\dots+h_n\log a_n| .
\end{equation}
\end{theorem}

\begin{corollary}\label{II:cor:main}
Let $a_1,\dots,a_n$ be as in Theorem~\ref{II:thm:main}.
\begin{enumerate}
\item $\rho(a)\le[\Q(a_1,\dots,a_n):\Q]$; hence \eqref{II:eq:main} holds with the exponent
$-(n-1+[\Q(a_1,\dots,a_n):\Q])-\varepsilon$.
\item If the $a_i$ are rational, the exponent is $n$; this is \cite[Corollary~1.2]{C}.
\item If $n=2$, the exponent is $\rho(a)+1$; this is also Corollary~\ref{I:cor:main}(5),
proved in Part~I by a separate argument.
\item If the $a_i$ lie in a subfield $K\subset\C$ of degree $d$, the exponent is at most
$n-1+d$; if they lie in one quadratic field, it is at most $n+1$.
\item If every $a_i$ is rational or a root of $x^2=tx-e$ with $t\in\Q$ and $e=\pm1$ (for
instance a unit of a real quadratic field), then $\rho(a)\le n$ and the exponent is at most
$2n-1$.
\end{enumerate}
\end{corollary}

Parts (2) and (4) are part~(1) with $[\Q(a):\Q]=1$ and $[\Q(a):\Q]\le d$; part~(3) is
Theorem~\ref{II:thm:main} for $n=2$. The rank can be much smaller than the degree:

\begin{example}\label{II:ex:quad}
Let $p_1,\dots,p_r$ be positive rational numbers without a power relation, and let
$u_1,\dots,u_s>1$ be units of $s$ pairwise different real quadratic fields. Then
$(p_1,\dots,p_r,u_1,\dots,u_s)$ has no power relation, and its rank is $s+1$ if $r\ge1$ and
$s$ if $r=0$, while the field it generates has degree $2^{s'}\ge s+1$ (it is the compositum
of the $s$ quadratic fields). So the exponent in \eqref{II:eq:main} is
$r+2s$ if $r\ge1$, and $2s-1$ if $r=0$. For instance, with
$\varphi=\frac{1+\sqrt5}2$,
\[
 |h_1\log2+h_2\log(1+\sqrt2)+h_3\log\varphi|\ \ge\ C\,H(h)^{-5-\varepsilon},
\]
while the degree gives the exponent $6$. This instance also follows from
Corollary~\ref{II:cor:main}(5); for $r\ge2$ the exponent $r+2s$ is smaller than the $2n-1$ of
part~(5).
\end{example}

The proof of the example is in Section~\ref{II:sec:cor}; the example (the absence of a power
relation and the computation of the rank) is not formalized.

\subsection{Two logarithms}

For real $x$ let
$\mu(x)=\sup\{\nu>0:\ 0<|x-p/q|<q^{-\nu}\text{ for infinitely many reduced }p/q\in\Q\}$.
For irrational $x$ it is invariant under $x\mapsto rx$, $r\in\Q^\times$, and under
$x\mapsto1/x$.

Let $a,b>0$ be real algebraic numbers, multiplicatively independent ($a^ub^v=1$ with
$u,v\in\Z$ only for $u=v=0$), so that $\xi=\log a/\log b$ is irrational. Let $\Pi(a,b)$ be
the group generated by the conjugate pairs $(\sigma b,\sigma a)$, so that
$\rk\Pi(a,b)=\rho(a,b)$ in the notation above.

\begin{theorem}\label{I:thm:main}
$2\le\mu(\log a/\log b)\le\rk\Pi(a,b)+2$.
\end{theorem}

\begin{corollary}\label{I:cor:main}
Let $a,b>0$ be multiplicatively independent real algebraic numbers.
\begin{enumerate}
\item $\mu(\log a/\log b)\le[\Q(a,b):\Q]+2$.
\item If $a,b$ lie in one quadratic field (one of them may be rational), then
$\mu(\log a/\log b)\le4$, and $\le3$ if both have norm $\pm1$.
\item If each of $a,b$ is rational or a root of $x^2=tx-\epsilon$ with $t\in\Q$,
$\epsilon=\pm1$ (for instance a quadratic unit), then $\mu(\log a/\log b)\le4$. For
instance, with $\varphi=\frac{1+\sqrt5}2$,
\[
 \mu\Bigl(\frac{\log2}{\log\varphi}\Bigr)\le4,\quad
 \mu\Bigl(\frac{\log(1+\sqrt2)}{\log\varphi}\Bigr)
 =\mu\Bigl(\frac{\artanh(1/\sqrt2)}{\artanh(1/\sqrt5)}\Bigr)\le4,\quad
 \mu\Bigl(\frac{\log(2+\sqrt3)}{\log(1+\sqrt2)}\Bigr)\le4 .
\]
\item For $\alpha=2\cos(2\pi/7)$ and its conjugate $\alpha'=\alpha^2-2=2\cos(4\pi/7)$,
$\mu(\log\alpha/\log2)\le5$ and $\mu(\log\alpha/\log|\alpha'|)\le4$.
\item For every $\varepsilon>0$ there is $c>0$ with
\[
 |p\log a+q\log b|\ge c\max(1,|p|,|q|)^{-\rk\Pi(a,b)-1-\varepsilon}
 \qquad\text{for all }(p,q)\in\Z^2\setminus0,
\]
and the same holds with $[\Q(a,b):\Q]$ in place of $\rk\Pi(a,b)$.
\end{enumerate}
\end{corollary}

The proofs are in Section~\ref{I:sec:cor}. Part~(2) with norms $\pm1$ is \cite{L}; for
rational $a,b$, $\rk\Pi=1$ and Theorem~\ref{I:thm:main} is \cite[Theorem~11.1]{L}.

The upper bound of Theorem~\ref{I:thm:main} and Corollary~\ref{I:cor:main}(5) are
equivalent, and (5) is the case $n=2$ of Theorem~\ref{II:thm:main}; in Lean it is derived
from the two-logarithm proof and, separately, from the $n$-logarithm proof
(Section~\ref{sec:lean}). These are two assemblies of one method, not independent
confirmations: both rest on the interpolation determinant of \cite{L}. Part~I proves more: Theorem~\ref{I:thm:B} bounds $\mu(\xi)$ by
$1+(t+2)/(z-t)$ for any family of $Q^z$ centers that lie on $O(Q^t)$ translates of the
logarithmic curve at every embedding, and it depends only on the two-torus interpolation
of \cite{L}. The examples (3) and (4) are proved with explicit families, without the Galois
closure.

\subsection{Context and method}

\paragraph{Relation to earlier work.}
By Dirichlet's theorem the exponent in \eqref{II:eq:main} cannot be smaller than $n-1$, for
any real numbers; for rational $a$ the bound $n$ comes from the theorem of \cite{C} for forms
with a constant term, where it is best possible. Baker's theory \cite{Ba,BW} gives effective
lower bounds whose exponents depend on the heights and degrees of the $a_i$, and effective
irrationality measures for $\log a/\log b$ of large size; see also \cite{W78}.
Theorem~\ref{II:thm:main} has an exponent that depends only on $n$ and $\rho(a)$, and an
ineffective constant. Its method is that of \cite{C} for rational numbers, which in turn is
the interpolation and determinant method of \cite{OAI} as developed in \cite{L}.

\cite[Section~12]{L} proves $\mu(\xi)\le3d/s$ when the logarithms match
($\lambda_{a,\sigma}=\xi\lambda_{b,\sigma}$) at $s$ of the $d$ embeddings of a field
containing $a,b$; the loss $d/s$ comes from the field norm of the interpolation
determinant, which is analytically small only at the matched embeddings. In the examples of
Corollary~\ref{I:cor:main}(3) at most half of the embeddings match, so \cite{L} gives $6$;
in part~(4) only the identity matches and \cite{L} gives $9$. As in the single-logarithm
case \cite{G}, we make every embedding good by enlarging the family of centers; the price
is the number of translated curves, or labels. Theorem~\ref{I:thm:main} is the analogue of
\cite[Theorem~1.1]{G}, $\mu(\log a)\le r(a)+2$: the conjugate ratios are replaced by the
conjugate pairs, and in Theorem~\ref{II:thm:main} by the conjugate tuples. Forms with a
constant term in logarithms of algebraic numbers are not treated here.

\paragraph{The idea for two logarithms.}
Let $\tau$ be an embedding and $(y,x)$ logarithms of $(\tau(b_c),\tau(a_c))$. The local
additive coordinates of \cite[Section~11]{L} are $\log(1+s)-r_i\log(1+t)$; in the global
coordinates $z=y+\log(1+t)$, $w=x+\log(1+s)$ this is
\[
 (w-\xi z-\kappa)+(\xi-r_i)\log(1+t),\qquad \kappa=x-\xi y .
\]
So every center lies exactly on the translated curve $w-\xi z=\kappa$, and the error
$(\xi-r_i)\log(1+t)$ is the same as in the rational case. The determinant at $\tau$ is
then estimated as in \cite{L}, except that the entire functions depend on $\kappa$: the
collision labels become pairs (translation, transverse index), and the three-halves
collision law loses the factor $N^{-1/2}$, $N$ the number of translations
(Theorem~\ref{I:thm:B}). The centers are products $\prod_lP_l^{c_l}$ of independent pairs
$P_l$ which at every embedding have a power of $(b,a)$ in their span. The translation of
$(b,a)$ is $\log a-\xi\log b=0$, so at every embedding the translation forgets one
direction of the box: $Q^z$ centers, about $Q^{z-1}$ translations, and
$1+(t+2)/(z-t)=z+2$. For a largest independent family of conjugate pairs,
$z=\rk\Pi(a,b)$.

\paragraph{The idea for $n$ logarithms.}
For three or more logarithms the additive coordinates are
$\log Y_n-\sum_jr_{i,j}\log Y_j$ with $m$ different slope vectors $r_i$, which approximate a
hyperplane rather than a point, and a center seems to carry $m$ translations, one for each
linear form. This is resolved as follows. The determinant of \cite{C} for $k+1=n$ logarithms
has its centers at the powers $\beta^J$ of one rational point of the torus
$\mathbb G_m^{k+1}$. We replace them by the products
\[
 b_c=\prod_{l}\beta_l^{c_l},\qquad c\in\{0,\dots,Q-1\}^z ,
\]
of a largest independent family $\beta_1,\dots,\beta_z$ of conjugate tuples of $a$. At an
embedding $\tau$ the conjugates $\tau(\beta_l)$ satisfy a relation
$\prod_l\tau(\beta_l)^{e_l}=a^{A}$; hence the logarithm vectors $Y_c$ of $\tau(b_c)$ can be
chosen as
\[
 Y_c=Y^0_{\mu(c)}+s_cL,\qquad L=(\log a_1,\dots,\log a_n),\ s_c\in\Z,
\]
with only $O(Q^{z-1})$ \emph{labels} $\mu(c)$: one direction of the box is absorbed by the
relation. The exact linear forms of \cite{C} vanish at $L$, so the entire functions of the
determinant at the center $c$ are $e^{h\cdot Y^0_{\mu(c)}}$ times the functions of
\cite{C} at the point $s_cL$ (Lemma~\ref{II:lem:shift}): all $m$ translations of a center
are determined by its label. The collision argument of \cite{C} then runs within the pairs
(label, grade), and loses the factor $N^{-1/n}$, $N$ the number of labels. With $K=Q^z$
centers, $N\asymp Q^{z-1}$ labels and the height condition $c^nB<1$, the parameters exist
exactly for $\nu>n+z$, which gives the exponent $n-1+z$ (Lemma~\ref{II:lem:params},
Remark~\ref{II:rem:where}).

\paragraph{Formal verification.}
Theorems~\ref{II:thm:main}, \ref{I:thm:main} and \ref{I:thm:B}, Corollary~\ref{II:cor:main},
Corollary~\ref{I:cor:main} (except the identities with $\artanh$ in (3)), and the lemmas and
propositions of the proofs are formally verified in Lean~4. Example~\ref{II:ex:quad}, the
numerical check of Remark~\ref{I:rem:num}, the remarks of Section~\ref{sec:remarks} and the
other points listed at the end of Section~\ref{sec:lean} are not. The formalization adds
41 files, 10572 lines, to the Lean library of \cite{C} (which contains that of \cite{L}):
21 files for Part~I and 20 for Part~II. It uses the interpolation theorems of \cite{L} and
\cite{C} as proved there, and edits no existing file. The main theorems depend only on the
axioms \lean{propext}, \lean{Classical.choice} and \lean{Quot.sound}. The proofs below follow
the formal proofs, and Lean names are given for the statements.

\paragraph{Organization and conventions.}
The two parts can be read independently, and each keeps the notation of its Lean files. In
both parts $E$ is a number field containing the numbers. In Part~I the two numbers are $b,a$,
with $\xi=\log a/\log b$. In Part~II $L$ is the vector of logarithms, $k\ge1$ is the number of \emph{slopes} and $k+1=n$
the number of logarithms and of torus coordinates; indices start at $0$, and the numbers are
$a_0,\dots,a_{k-1},g$, with logarithms $L_j=\log a_j$ $(j<k)$ and $G=L_k=\log g$.

\part{Two logarithms}\label{part:I}

\section{Inputs from \texorpdfstring{\cite{L}}{[L]}}\label{I:sec:inputs}

Fix $m\ge2$. Coordinates $Y,Z,X_1,\dots,X_m$ have weights $(w_0,w_0,w_1,\dots,w_m)$; jets
in $t,s,u_1,\dots,u_m$ have weights $(v,v,w_1/\theta,\dots,w_m/\theta)$. Write
$w\cdot\gamma=\sum_iw_i\gamma_i$, $w_*=\min_iw_i$. A polynomial $Y^hZ^kX^\gamma$ has
weight $w_0(h+k)+w\cdot\gamma$.

\begin{theorem}[{Two-torus interpolation, \cite[Lemma~11.2]{L}}]\label{I:thm:interp}
Let $K\ge1$, $0<\theta<1$, and $w_0,v,\sigma>0$ be rational with
\begin{equation}\label{I:eq:volume}
\begin{gathered}
 (1+3\sigma)^{m+2}K(w_0/v)^2\theta^m<1,\qquad
 (1+3\sigma)^{m+1}(w_0/v)\theta^{m-1}<1,\\
 (1+3\sigma)^{m-1}\theta^{m-1}<1,\qquad(1+\sigma)\theta<1,\qquad
 4(m+3)/\sigma<w_0/v .
\end{gathered}
\end{equation}
There are successive thresholds for rational weights $w_1,\dots,w_m$, each depending only
on these data and the earlier weights, such that for all weights above them there is
$R_0\ge1$ with the following property. For all pairwise distinct non-zero rationals
$r_1,\dots,r_m$ and all centers $(b_c,a_c)\in(\C^\times)^2$, $c=1,\dots,K$, such that for
every $i$ the characters $a_c^{q_i}/b_c^{p_i}$ ($r_i=p_i/q_i$ reduced, $q_i\ge1$) are
pairwise distinct, there is $H_0$ such that
\[
 P\longmapsto\Bigl[P\bigl(b_c(1+t),\,a_c(1+s),\,(\log(1+s)-r_i\log(1+t)+u_i)_i\bigr)
 \bmod F_V^{\ge H}\Bigr]_{c\le K}
\]
from polynomials of weight $\le H$ onto all jets of weight $<H$ at all centers is
surjective for every $H\ge H_0$ divisible by $R_0$.
\end{theorem}

In Lean this is the statement \lean{TwoLogInterpolation.UniformStatement}, proved in the
file \lean{TwoLogUniformInterpolation}; it is the statement displayed, with arbitrary centers
separated by characters. In \cite[Lemma~11.2]{L} it is written for the centers
$(3^j,2^j)$; its proof uses the centers only through the distinctness of the characters.

\begin{lemma}[{Truncation, \cite[Lemma~11.3]{L}}]\label{I:lem:shear}
Let $T_i$ be integers with $vT_i>w_i/\theta$ and $G_i(t)=\sum_{1\le l<T_i}(-1)^{l+1}t^l/l$.
For any fixed set of centers and any selection of as many columns as there are rows, the
coefficient matrices built with $G_i(s)-r_iG_i(t)+u_i$ and with
$\log(1+s)-r_i\log(1+t)+u_i$ have the same determinant.
\end{lemma}

\begin{lemma}[{Bivariate collision, \cite[Lemma~11.4]{L}}]\label{I:lem:coll}
For $n$ distinct pairs $(k,l)\in\Z_{\ge0}^2$, $\sum(k+l)\ge\frac{2\sqrt2}3n^{3/2}-n$.
\end{lemma}

We use $\log\lcm(1,\dots,n)\le\Lambda n$ with $\Lambda=\log4+4$, and the weighted
lattice-point asymptotics of \cite[Section~11.2]{L}.

\section{Two-torus center supplies}\label{I:sec:supply}

\begin{definition}\label{I:def:supply}
Let $\xi\in\R$, $E$ a number field and $z,t\ge0$ integers. A \emph{two-torus center supply
of type $(z,t)$} for $\xi$ in $E$ consists of integers $D,c_N\ge1$, reals
$c_E,c_R\ge0$, $c_S$, and, for every integer $Q\ge1$, centers
$(b_c,a_c)\in(E^\times)^2$ indexed by a finite set $\mathcal C_Q$ and an integer $e\ge0$, such that:
\begin{enumerate}
\item[(T1)] $|\mathcal C_Q|=Q^z$;
\item[(T2)] $e\le c_EQ$, and $D^eb_c$, $D^ea_c$ are algebraic integers;
\item[(T3)] for every embedding $\tau:E\to\C$ there are $y_\tau(c),x_\tau(c)\in\C$ with
$e^{y_\tau(c)}=\tau(b_c)$, $e^{x_\tau(c)}=\tau(a_c)$, and
$|y_\tau(c)|,|x_\tau(c)|\le c_RQ$;
\item[(T4)] for every embedding $\tau$, the \emph{translations} $\kappa_\tau(c)=x_\tau(c)-\xi y_\tau(c)$ take at most
$c_NQ^t$ values (each of modulus $\le(1+|\xi|)c_RQ$ by (T3));
\item[(T5)] for every reduced $p/q\ne0$ with $q>c_SQ$, the map
$c\mapsto a_c^qb_c^{-p}$ is injective on $\mathcal C_Q$.
\end{enumerate}
\end{definition}

\begin{theorem}[Theorem B]\label{I:thm:B}
Let $\xi$ be irrational with a two-torus center supply of type $(z,t)$, $s=z-t$,
$1\le s\le t+2$. Then for every $\nu>1+(t+2)/s$ there are only finitely many reduced $p/q$
with $q\ge2$ and $0<|\xi-p/q|<q^{-\nu}$. In particular
$2\le\mu(\xi)\le1+\dfrac{t+2}{z-t}$.
\end{theorem}

The condition $s\le t+2$ says that the bound is at least $2$. Sections~\ref{I:sec:construction}--\ref{I:sec:end}
prove Theorem~\ref{I:thm:B}. In this section we construct the supplies.

\paragraph{Conjugate pairs.}
In the proof of Theorem~\ref{I:thm:main}, $E$ is a Galois number field with an embedding
$\iota:E\to\C$ and elements $b_E,a_E$ with $\iota(b_E)=b$, $\iota(a_E)=a$. The conjugate
pairs are $g(b_E,a_E)=(gb_E,ga_E)$, $g\in\Gal(E/\Q)$. A family $g_1,\dots,g_z$ is
\emph{independent} if $\prod_lg_l(b_E,a_E)^{n_l}=(1,1)$ implies $n=0$; the largest $z$ of
an independent family is $\rk\Pi(a,b)$.

\begin{lemma}[Relation]\label{I:lem:relation}
Let $g_1,\dots,g_z$ be an independent family of the largest size and $P_l=g_l(b_E,a_E)$.
For every $g\in\Gal(E/\Q)$ there are $A\in\Z\setminus0$ and $n\in\Z^z$ with
$\prod_lP_l^{n_l}=g(b_E,a_E)^{A}$.
\end{lemma}

\begin{proof}
The family $g_1,\dots,g_z,g$ is not independent; in a non-trivial relation the exponent of
$g(b_E,a_E)$ is non-zero since $g_1,\dots,g_z$ are independent, and its negative is $A$.
\end{proof}

\begin{proposition}[The supply]\label{I:prop:supply}
Let $b,a>0$, $\xi=\log a/\log b$ irrational, and let $E$ be a number field with
pairs $P_l=(\beta_l,\alpha_l)\in(E^\times)^2$, $l=1,\dots,z$, such that
\begin{enumerate}
\item[(i)] the $P_l$ are independent: $\prod_lP_l^{n_l}=(1,1)$ implies $n=0$;
\item[(ii)] $D\beta_l$, $D\alpha_l$ are algebraic integers for some integer $D\ge1$;
\item[(iii)] for every embedding $\sigma:E\to\C$ there are $A_\sigma\in\Z$ and
$n\in\Z^z\setminus0$ with $\prod_l\sigma(\beta_l)^{n_l}=b^{A_\sigma}$ and
$\prod_l\sigma(\alpha_l)^{n_l}=a^{A_\sigma}$.
\end{enumerate}
For $c\in\{0,\dots,Q-1\}^z$ put $(b_c,a_c)=\prod_lP_l^{c_l}$. This is a two-torus center
supply of type $(z,z-1)$ for $\xi$ in $E$; hence $\mu(\xi)\le z+2$.
\end{proposition}

We call such pairs a \emph{pair family of size $z$} for $(b,a)$.

\begin{proof}
Fix $\sigma$, logarithms $B_{\sigma,l}$ of $\sigma(\beta_l)$ and $A_{\sigma,l}$ of
$\sigma(\alpha_l)$, and the relation (iii), after a change of sign of $(n,A_\sigma)$ with
$n_{l_0}>0$ for some $l_0$. There are integers $k_b,k_a$ with
\[
 \sum_ln_lB_{\sigma,l}=A_\sigma\log b+2\pi ik_b,\qquad
 \sum_ln_lA_{\sigma,l}=A_\sigma\log a+2\pi ik_a .
\]
For $c$ in the box let $j=\lfloor c_{l_0}/n_{l_0}\rfloor$ and
\[
 y_\sigma(c)=\sum_lc_lB_{\sigma,l}-2\pi ijk_b,\qquad
 x_\sigma(c)=\sum_lc_lA_{\sigma,l}-2\pi ijk_a .
\]

\emph{(T2).} With $e=zQ$, $D^eb_c=D^{e-\sum c_l}\prod_l(D\beta_l)^{c_l}$ is integral;
likewise $D^ea_c$.

\emph{(T3).} $e^{y_\sigma(c)}=\prod_l\sigma(\beta_l)^{c_l}=\sigma(b_c)$, likewise for $x$;
$0\le j<Q$ gives $|y_\sigma(c)|,|x_\sigma(c)|\le c_RQ$ with $c_R$ depending on the
logarithms and on $k_a,k_b$ for all $\sigma$.

\emph{(T4).} Put $c^*=c-jn$. Since $\log a=\xi\log b$,
\begin{align*}
 \kappa_\sigma(c)&=\sum_lc_l(A_{\sigma,l}-\xi B_{\sigma,l})-2\pi ij(k_a-\xi k_b)
 =\sum_lc^*_l(A_{\sigma,l}-\xi B_{\sigma,l})+jA_\sigma(\log a-\xi\log b)\\
 &=\sum_lc^*_l(A_{\sigma,l}-\xi B_{\sigma,l}),
\end{align*}
so $\kappa_\sigma(c)$ depends only on $c^*$. Let $N_b$ bound all $|n_l|$, for all
$\sigma$. Then $0\le c^*_{l_0}<n_{l_0}\le N_b$ and $|c^*_l|\le Q+QN_b$, hence
$\#\{\kappa_\sigma(c)\}\le N_b(2N_b+3)^{z-1}Q^{z-1}$ (Lean uses the cruder
$c_N=(N_b+1)(2N_b+3)^{z-1}$).

\emph{(T5).} Let $G\subset E^\times$ be the group generated by all $\beta_l,\alpha_l$.
Modulo its torsion it is free of finite rank; fix coordinates $\phi:G\to\Z^g$ with kernel
the torsion, and let $c_S=z\max_l(\|\phi(\beta_l)\|_\infty+\|\phi(\alpha_l)\|_\infty)$.
Suppose $a_c^qb_c^{-p}=a_{c'}^qb_{c'}^{-p}$ with $\gcd(p,q)=1$, $q>c_SQ$, and put
$\Delta=c-c'$, $\mathsf A=\prod_l\alpha_l^{\Delta_l}$, $\mathsf B=\prod_l\beta_l^{\Delta_l}$.
Then $\mathsf A^q=\mathsf B^p$. With $xp+yq=1$ and $u=\mathsf A^x\mathsf B^y\in G$ we get
$u^p=\mathsf A$ and $u^q=\mathsf B$. So $\phi(\mathsf B)=q\phi(u)$, while
$\|\phi(\mathsf B)\|_\infty<c_SQ<q$; hence $\phi(u)=0$, $u$ is torsion, and so are
$\mathsf A$ and $\mathsf B$. Then $\prod_lP_l^{N\Delta_l}=(1,1)$ for some $N\ge1$, and
$\Delta=0$ by (i).

With $t=z-1$, $s=1\le t+2$, Theorem~\ref{I:thm:B} gives $\mu(\xi)\le1+(z+1)=z+2$.
\end{proof}

\begin{proof}[Proof of Theorem~\ref{I:thm:main}]
The maps $g\mapsto\iota\circ g$ from $\Gal(E/\Q)$ to the embeddings of $E$ are
injective, and both sets have $[E:\Q]$ elements; so every embedding is $\sigma=\iota\circ
g$. Let $g_1,\dots,g_z$ be an independent family of the largest size, $P_l=g_l(b_E,a_E)$,
and apply Lemma~\ref{I:lem:relation} to $g^{-1}$: $\prod_lP_l^{n_l}=g^{-1}(b_E,a_E)^A$, and
applying $\sigma=\iota\circ g$ gives (iii) with $A_\sigma=A$. Here $n\ne0$: otherwise
$b^A=1$ with $A\ne0$, impossible as $b\ne1$. With $D$ a common denominator of $b_E,a_E$,
the conjugates $g_lb_E,g_la_E$ satisfy (ii). Proposition~\ref{I:prop:supply} gives
$\mu\le z+2$; the lower bound is Dirichlet's theorem.
\end{proof}

\begin{lemma}\label{I:lem:cubic}
Let $\alpha=\alpha_1=2\cos(2\pi/7)$, $\alpha_2=\alpha^2-2$, $\alpha_3=1-\alpha-\alpha^2$,
the roots of $x^3+x^2-2x-1$, and $E=\Q(\theta)$, $\theta^3+\theta^2-2\theta-1=0$, with
$\theta_1=\theta$, $\theta_2=\theta^2-2$, $\theta_3=1-\theta-\theta^2$. Then
$(2,\theta_k)$, $k=1,2,3$, is a pair family of size $3$ for $(2,\alpha)$, and
$(-\theta_2,\theta_1)$, $(-\theta_3,\theta_2)$ is a pair family of size $2$ for
$(|\alpha_2|,\alpha)=(-\alpha_2,\alpha_1)$.
\end{lemma}

\begin{proof}
The cubic is irreducible (its only possible rational roots $\pm1$ fail), and the three
embeddings send $(\theta_1,\theta_2,\theta_3)$ to the cyclic shifts of
$(\alpha_1,\alpha_2,\alpha_3)$, since $\alpha_2=2\cos(4\pi/7)$, $\alpha_3=2\cos(6\pi/7)$
and $\cos4t=\cos3t$ at $t=2\pi/7$. The $\alpha_k$ are units with $\alpha_1\alpha_2\alpha_3=1$.
Put $\ell_k=\log|\alpha_k|$, so $\ell_1+\ell_2+\ell_3=0$. If $\alpha_1^x\alpha_2^y=\pm1$,
then $x\ell_1+y\ell_2=0$ and, applying the cyclic embedding, $x\ell_2+y\ell_3=0$; the
determinant is $-(\ell_1^2+\ell_1\ell_2+\ell_2^2)\ne0$ as $\ell_1\ne0$, so $x=y=0$.
\emph{First family.} If $\prod_k(2,\theta_k)^{n_k}=(1,1)$ then $n_1+n_2+n_3=0$ and
$\theta_1^{n_1-n_3}\theta_2^{n_2-n_3}=1$, so all $n_k$ are equal and hence $0$. At the
embedding sending $\theta_k$ to $\alpha_1$ the relation is $n=e_k$, $A=1$.
\emph{Second family.} Independence follows from the first step. At the identity, $n=(1,0)$;
at $\theta_1\mapsto\alpha_3$, $n=(0,1)$; at $\theta_1\mapsto\alpha_2$ the pairs go to
$(-\alpha_3,\alpha_2)$, $(-\alpha_1,\alpha_3)$, whose product is
$(\alpha_2^{-1},\alpha_1^{-1})$, so $n=(-2,-2)$, $A=2$.
\end{proof}

\begin{remark}[Numerical check]\label{I:rem:num}
For $(\varphi^2,(1+\sqrt2)^2)$ in $\Q(\sqrt2,\sqrt5)$ the four embeddings send the pair to
the four sign patterns $(b^{\pm1},a^{\pm1})$, and the grid $(b^{j_1},a^{j_2})$ has exactly
$2Q-1$ translations at each of them ($Q=5,10,20,40$). For $(\alpha^2,2)$ the box of
centers $(2^{c_0},a_1^{c_1}a_2^{c_2})$, $a_k=\alpha_k^2$, gives $496,721,496$ translations
at the three embeddings for $Q=16$ ($4096$ centers); for the conjugate units $(a_1,a_2)$
the centers of Proposition~\ref{I:prop:supply} give $Q,Q,2Q-1$ translations for $Q^2$
centers ($Q=10,20,40$).
\end{remark}

\subsection{Proof of Corollary~\ref{I:cor:main}}\label{I:sec:cor}

\emph{(1)} Let first $E$ be Galois as above and $F=\Q(b_E,a_E)$. In an independent family
the pairs $g_l(b_E,a_E)$ are distinct (two equal pairs give the relation
$P_iP_j^{-1}=(1,1)$), so the restrictions $g_l|_F:F\to E$ are distinct $\Q$-algebra maps.
By Dedekind--Artin they are linearly independent over $E$ in $\operatorname{Hom}_\Q(F,E)$,
whose dimension over $E$ is $[F:\Q]$. Hence $\rk\Pi(a,b)\le[F:\Q]$. For real $a,b$ take
for $E$ the splitting field of the product of the minimal polynomials of $a$ and $b$,
embedded in $\C$ by some $\iota$; the polynomial splits in $E$, so $a,b$ have preimages
$a_E,b_E$, and $\iota$ maps $\Q(b_E,a_E)$ isomorphically onto $\Q(b,a)$.

\emph{(2)} Let $E$ be the quadratic field, $\iota$ its embedding with $\iota(b_E,a_E)=(b,a)$
and $\sigma$ the other one. For $x\in E$, $\iota(x)\sigma(x)=N(x)$, so
$x'=N(x)/x\in E$ satisfies $\sigma(x')=\iota(x)$. If $(b_E,a_E)$ and $(b'_E,a'_E)$ are
independent, they form a pair family of size $2$ (relations $n=(1,0)$ at $\iota$ and
$n=(0,1)$ at $\sigma$, $A=1$), and $\mu\le4$. Otherwise there is a relation
$(b_E,a_E)^{n_0}(b'_E,a'_E)^{n_1}=(1,1)$ with $n_0\ne0$ (if $n_0=0$, then $n_1\ne0$ and
$b=\sigma(b'_E)$ would be a root of unity, while $b>0$, $b\ne1$); applying $\sigma$ gives
$\sigma(b_E,a_E)^{n_0}=(b,a)^{-n_1}$, so $(b_E,a_E)$ alone is a pair family of size $1$,
and $\mu\le3$. If both norms are $\pm1$ then $(b'_E,a'_E)=(\pm b_E^{-1},\pm a_E^{-1})$ and
the second case occurs.

\emph{(3)} Let $F=\Q(a,b)\subset\R$. For every embedding $\sigma$ of $F$ and $x\in\{a,b\}$,
$\sigma(x)^2=x^{\pm2}$: $\sigma(x)=x$ if $x$ is rational, and otherwise $\sigma(x)$ is a
root of $x^2=tx-\epsilon$, whose roots are $x$ and $\epsilon/x$. So $(b,1)$, $(1,a)$ is a
pair family of size $2$ (independent as $a,b$ are; at $\sigma$ the relation is
$n=(\pm2,\pm2)$, $A=2$), and $\mu\le4$. In the examples, $2^u\varphi^v=1$ gives, after
conjugation, $2^u\bar\varphi^v=1$ and $4^u(-1)^v=1$, so $u=v=0$; two units $w_1,w_2>1$ of
different real quadratic fields are multiplicatively independent, since every power is
$w^v=A_v+B_vw$, $w'^v=A_v+B_vw'$ with integers $A_v,B_v$, $B_v\ne0$ for $v\ne0$, and
$w_1^u=w_2^v$ would make $w_1-w_1'$ or $(w_1-w_1')(w_2-w_2')$ rational. Finally
$\log\varphi=\artanh(1/\sqrt5)$, $\log(1+\sqrt2)=\artanh(1/\sqrt2)$.

\emph{(4)} Lemma~\ref{I:lem:cubic} and Proposition~\ref{I:prop:supply} with $z=3$ and $z=2$;
the multiplicative independence of $\alpha$ and $|\alpha_2|$ is shown in the proof of
Lemma~\ref{I:lem:cubic}, and that of $\alpha$ and $2$ follows from the norms: $\alpha^x2^y=1$
gives $8^y=N(\alpha)^{-x}=1$, so $y=0$, and then $x=0$ since $\ell_1\ne0$.

\emph{(5)} The proofs above give the finiteness statement of Theorem~\ref{I:thm:B}: with
$\kappa=\rk\Pi(a,b)+2$ (or $[\Q(a,b):\Q]+2$), every $\nu>\kappa$ has only finitely many good
approximations. Put $\nu=\kappa+\varepsilon$, $\xi=\log a/\log b$, and for $p\ne0$ let
$r=-q/p$ in lowest terms, $d$ its denominator, so $d\le|p|$. Let $N\ge2$ exceed the
denominators of the finitely many good approximations. If $d\ge N$, then $r$ is not a good
approximation and
$|\xi-r|\ge d^{-\nu}\ge|p|^{-\nu}$. If $d<N$ then $|d\xi-dr|\ge\delta$ for a $\delta>0$
depending on $N$ ($\xi$ is irrational), so $|p\xi+q|=|p||\xi-r|\ge\delta|p|/d\ge\delta$.
In both cases $|p\xi+q|\ge\min(1,\delta)|p|^{1-\nu}$. As $p\log a+q\log b=\log b\,(p\xi+q)$,
and $|q\log b|\ge|\log b|$ for $p=0$, this gives (5) with exponent $1-\nu=-(\kappa-1+\varepsilon)$. \qed

\section{The construction}\label{I:sec:construction}

Let $\xi$ have a supply of type $(z,t)$ in $E$ with denominator $D$, $s=z-t$,
$1\le s\le t+2$. Fix $\nu$ with $\nu s>z+2$, i.e.\ $\nu>1+(t+2)/s$; then $\nu>2$. Suppose,
to reach a contradiction, that there are infinitely many reduced $p/q$, $q\ge2$, with
$0<|\xi-p/q|<q^{-\nu}$; then there are such $p/q$ with arbitrarily large $q$.

\paragraph{Fixed parameters.}
Choose rationals $\theta,A,B,c$ as in Lemma~\ref{I:lem:params}:
\begin{equation}\label{I:eq:P}
 0<\theta<A<1,\quad\nu(A-\theta)>1-\theta,\quad
 Ac^t<B<\min\bigl(c^{-2},\,c^z\theta\bigr),\quad c^z\theta<1,\quad c>1 .
\end{equation}
Choose $0<\eta<1$ with $g:=\nu(A(1-\eta)-\theta)-(1-\theta)>0$, and put $\epsilon=g/2$ and
$F_0=4/\theta>2/\theta$.

\paragraph{Dimension and centers.}
For $m=2n$, fixed in Section~\ref{I:sec:end}, put $Q=\lfloor c^{2n}\rfloor$ and take the
$K=Q^z$ centers $(b_c,a_c)$ supplied for this $Q$. Put
\begin{equation}\label{I:eq:m}
 S=\lfloor\sqrt K\rfloor+1,\quad w_0=B^{-n},\quad v=2S\theta^nB^{-n},\quad
 \sigma=\frac{1-\theta}{100(m+3)},\quad N=c_NQ^t,\quad R=2c_RQ+2 .
\end{equation}
Then $w_0,v,\sigma$ are rational and $K(w_0/v)^2\theta^m=K/(4S^2)\le\frac14$.

\paragraph{Approximations.}
With $m$ fixed, choose distinct non-zero reduced fractions $r_i=p_i/q_i$, $q_i\ge2$, with
$|\xi-r_i|<q_i^{-\nu}$, successively for $i=1,\dots,m$, so that
$w_i=\lceil\log q_i\rceil$ exceeds the $i$-th threshold of Theorem~\ref{I:thm:interp}, a
bound $X+1$ fixed in Section~\ref{I:sec:end}, and $\lceil\log(\lceil c_SQ\rceil+1)\rceil$;
the last gives $q_i>c_SQ$. Put
\[
 T_i=\lceil F_0w_i/v\rceil,\quad G_i(t)=\sum_{1\le l<T_i}\frac{(-1)^{l+1}t^l}l,\quad
 L_i=\lcm(1,\dots,T_i-1),\quad \epsilon_i(t)=(\xi-r_i)\log(1+t).
\]
Then $vT_i\ge F_0w_i>w_i/\theta$.

\paragraph{Matrix and minor.}
Columns: $P=Y^hZ^kX^\gamma$ with $w_0(h+k)+w\cdot\gamma\le H$. Rows: $(c,\lambda,\rho,\beta)$
with $c$ a center and $v(\lambda+\rho)+w\cdot\beta/\theta<H$. Entry:
\begin{equation}\label{I:eq:entry}
 [t^\lambda s^\rho u^\beta]\ P\bigl(b_c(1+t),\,a_c(1+s),\,(G_i(s)-r_iG_i(t)+u_i)_i\bigr)\in E .
\end{equation}
Fix an embedding $\sigma_0$ of $E$. By (T5) and $q_i>c_SQ$ the characters of the centers
$(\sigma_0b_c,\sigma_0a_c)$ are distinct, and \eqref{I:eq:volume} holds
(Section~\ref{I:sec:end}). Theorem~\ref{I:thm:interp} and Lemma~\ref{I:lem:shear} give, for
arbitrarily large $H$, a selection of columns such that the square minor
$\Delta=\Delta_H\in E$ using all $M=M_H$ rows has $\sigma_0(\Delta)\ne0$; hence $\Delta\ne0$
and $\tau(\Delta)\ne0$ for every $\tau$. Put $\bar b=\frac1{MH}\sum_{\rm rows}w\cdot\beta\in
[0,\theta)$. As $H\to\infty$ with everything else fixed,
\begin{equation}\label{I:eq:counts}
 M\sim\frac{K\theta^mH^{m+2}}{(m+2)!\,v^2\prod_iw_i},\qquad
 D_A(H):=\#\{\delta\in\Z_{\ge0}^m:w\cdot\delta\le AH\}\sim\frac{A^mH^m}{m!\prod_iw_i}.
\end{equation}

\section{Arithmetic lower bound}\label{I:sec:arith}

\begin{lemma}\label{I:lem:arith}
\[
 \frac1{[E:\Q]}\sum_\tau\frac{\log|\tau\Delta|}{MH}\ \ge\ -(1-\bar b)-E_{\rm ar},\qquad
 E_{\rm ar}=\frac{\Lambda F_0}{v}+\frac{2e\log D}{w_0},
\]
the sum over all embeddings $\tau$ of $E$.
\end{lemma}

\begin{proof}
A non-zero term of the determinant expansion picks entries with $\beta\le\gamma$, so
$J_i=\sum_{\rm col}\gamma_i-\sum_{\rm rows}\beta_i\ge0$ \cite[Section~11.2]{L}. Scale
column $(h,k,\gamma)$ by $(D^e)^{2\lfloor H/w_0\rfloor}\prod_i(q_iL_i)^{\gamma_i}$ and row
$\beta$ by $\prod_i(q_iL_i)^{-\beta_i}$. A non-zero entry equals
$b_c^ha_c^k\binom\gamma\beta[t^\lambda s^\rho](1+t)^h(1+s)^k
\prod_i(G_i(s)-r_iG_i(t))^{\gamma_i-\beta_i}$; by (T2) and $h,k\le H/w_0$,
$(D^e)^{2\lfloor H/w_0\rfloor}b_c^ha_c^k$ is integral, and
$q_iL_i(G_i(s)-r_iG_i(t))=q_iL_iG_i(s)-p_iL_iG_i(t)\in\Z[s,t]$. Hence $C_H\Delta$ is a
non-zero algebraic integer, $C_H=D^{2eM\lfloor H/w_0\rfloor}\prod_i(q_iL_i)^{J_i}$, and its
norm gives $\sum_\tau\log|\tau\Delta|\ge-[E:\Q]\log C_H$. Finally $\sum_iw_iJ_i\le
MH(1-\bar b)$, $\log q_i\le w_i$ and $\log L_i\le\Lambda(T_i-1)<\Lambda F_0w_i/v$.
\end{proof}

\section{Analytic upper bound at every embedding}\label{I:sec:analytic}

For a column $P$, $\kappa\in\C$ and $\delta\in\Z_{\ge0}^m$ define the entire function
\[
 f^\kappa_{\delta,P}(z,w)=[u^\delta]\,P\bigl(e^z,e^w,(w-\xi z-\kappa+u_i)_i\bigr).
\]

\begin{lemma}[Row identity]\label{I:lem:row}
Let $\tau$ be an embedding, $c$ a center, $y=y_\tau(c)$, $x=x_\tau(c)$,
$\kappa=\kappa_\tau(c)$. For every polynomial $P$,
\begin{multline*}
 [u^\beta]\,P\bigl(\tau(b_c)(1+t),\tau(a_c)(1+s),(\log(1+s)-r_i\log(1+t)+u_i)_i\bigr)\\
 =\sum_{\delta\ge\beta}\binom\delta\beta\epsilon(t)^{\delta-\beta}
 f^\kappa_{\delta,P}\bigl(y+\log(1+t),\,x+\log(1+s)\bigr).
\end{multline*}
\end{lemma}

\begin{proof}
With $z=y+\log(1+t)$, $w=x+\log(1+s)$, (T3) gives $e^z=\tau(b_c)(1+t)$,
$e^w=\tau(a_c)(1+s)$, and
$w-\xi z-\kappa=\log(1+s)-\xi\log(1+t)$; adding $\epsilon_i(t)$ gives
$\log(1+s)-r_i\log(1+t)$. Expand in $u$ by the binomial formula.
\end{proof}

The $r_i$ are rational, so $\tau$ fixes them, and by Lemma~\ref{I:lem:shear} the matrix of
$\tau\Delta$ may be computed with full logarithms; its rows are exactly the left sides of
Lemma~\ref{I:lem:row}. \emph{Every embedding is good}: only the label $\kappa$ depends on
$\tau$, and by (T4) it takes at most $N$ values.

\begin{proposition}\label{I:prop:upper}
Let $c_2=\frac{\log2}8$. For every embedding $\tau$ and all $H$ as in
Section~\ref{I:sec:construction},
\[
 \frac{\log|\tau\Delta|}{MH}\le E_{\rm an}+o(1)
 +\max\Bigl\{-c_2\,\eta^{3/2}\sqrt{\frac{M}{H^2\,N\,D_A(H)}},\ -\nu\bigl(A(1-\eta)-\bar b\bigr)
 \Bigr\},
\]
where the term $o(1)$, as $H\to\infty$, does not depend on $\tau$, and
\[
 E_{\rm an}=\frac R{w_0}+\frac{\log2}v+\frac{\log B_R+\nu+\log2}{w_*},\qquad
 B_R=4R(1+|\xi|) .
\]
\end{proposition}

\begin{proof}
\emph{Expansion.} By Lemma~\ref{I:lem:row}, extracting $[t^\lambda s^\rho]$, the entry in row
$(c,\lambda,\rho,\beta)$ and column $P$ is
\[
 \sum_{\delta\ge\beta}\sum_{\ell\le\lambda}\binom\delta\beta[t^\ell]\epsilon(t)^{\delta-\beta}
 \cdot[t^{\lambda-\ell}s^\rho]\,f^{\kappa}_{\delta,P}\bigl(y+\log(1+t),x+\log(1+s)\bigr),
\]
with $\kappa=\kappa_\tau(c)$; only $w\cdot\delta\le H$ occurs. Expand
$f^\kappa_{\delta,P}=\sum_{k,l}F^\kappa_{\delta,k,l}(P)(z/R)^k(w/R)^l$. By multilinearity
$\tau\Delta$ is a sum, over choices of $(\delta,\ell,k,l)$ for every row, of products of
the row scalars and the determinants $\det[F^{\kappa_r}_{\delta_r,k_r,l_r}(P)]_{r,P}$. Such
a determinant vanishes unless the quadruples $(\kappa_r,\delta_r,k_r,l_r)$ are pairwise
distinct: rows with the same quadruple are equal. So within each \emph{label}
$(\kappa,\delta)$ the pairs $(k,l)$ are distinct.

\emph{Sizes.} By (T3), $|y|+\frac34\le R/2$ and $|x|+\frac34\le R/2$, and
$|\kappa|\le(1+|\xi|)c_RQ\le R(1+|\xi|)$. The points $(y+\log(1+t),x+\log(1+s))$ with
$|t|=|s|=\frac12$ lie in the polydisc of radius $R/2$, since $|\log(1+t)|\le\log2<\frac34$.
On the polydisc of radius $R$, $|w-\xi z-\kappa+u_i|\le B_R$ for $|u_i|=1$; with
$h+k\le H/w_0$ and $|\gamma|\le H/w_*$, Cauchy's estimate gives
$|F^\kappa_{\delta,k,l}(P)|\le\exp(H(R/w_0+\log B_R/w_*))$. The row test
$[t^{\lambda-\ell}s^\rho]((y+\log(1+t))/R)^k((x+\log(1+s))/R)^l$ is at most
$2^{\lambda-\ell+\rho}2^{-k-l}$. The scalar satisfies, since $|\xi-r_i|\le q_i^{-\nu}\le
e^{-\nu(w_i-1)}$ and $\binom\delta\beta\le2^{|\delta|}$,
\[
 \Bigl|\binom\delta\beta[t^\ell]\epsilon(t)^{\delta-\beta}\Bigr|\le
 2^\ell\exp\Bigl(-\nu\,w\cdot(\delta-\beta)+\frac{(\nu+\log2)H}{w_*}\Bigr),
\]
and $2^{\lambda+\rho}\le e^{H\log2/v}$. The number of choices of $(\delta,\ell)$ per row,
the factor $M!$ and the geometric sums over $(k,l)$ (half of $2^{-k-l}$ is reserved)
contribute $o(MH)$; these give $E_{\rm an}+o(1)$.

\emph{Dichotomy.} Fix one choice of $\delta_r$ for all rows. If at least $\eta M$ rows have
$w\cdot\delta_r\le AH$, their labels $(\kappa,\delta)$ take at most $N\,D_A(H)$ values; if
$n_\ell$ rows carry label $\ell$, Lemma~\ref{I:lem:coll} and the reserved half of
$2^{-k-l}$ save $\exp(-c\sum_\ell(\frac{2\sqrt2}3n_\ell^{3/2}-n_\ell))$ with an absolute
$c>0$, and by H\"older's inequality
\[
 \sum_\ell n_\ell^{3/2}\ \ge\ \frac{(\eta M)^{3/2}}{\sqrt{N\,D_A(H)}} ,
\]
which gives the first term of the maximum (the linear term $\sum n_\ell\le M$ is $o(MH)$).
Otherwise more than $(1-\eta)M$ rows have $w\cdot\delta_r>AH$, and
$\sum_r w\cdot(\delta_r-\beta_r)\ge MH(A(1-\eta)-\bar b)$, which gives the second term.
In each alternative the other factor is at most $1$. Summing over the choices gives the
claim. This is \cite[Section~11.3]{L} with the label $\delta$ replaced by $(\kappa,\delta)$;
the library's constant is $c_2=\frac{\log2}8$.
\end{proof}

By \eqref{I:eq:counts}, as $H\to\infty$,
\begin{equation}\label{I:eq:Lprime}
 c_2\eta^{3/2}\sqrt{\frac M{H^2ND_A(H)}}\ \longrightarrow\ L'=
 c_2\eta^{3/2}\sqrt{\frac{K\theta^m}{(m+1)(m+2)\,v^2A^mN}} .
\end{equation}

\section{Parameters and the end of the proof}\label{I:sec:end}

\begin{lemma}\label{I:lem:params}
Let $z>t\ge0$ be integers and $\nu>1+(t+2)/(z-t)$. There are rationals $\theta,A,B,c$
with \eqref{I:eq:P}.
\end{lemma}

\begin{proof}
Put $s=z-t$. The interval $\frac{t+2}{z+2}<b<1-\frac1\nu$ is non-empty, since
$1-\frac1\nu>1-\frac{s}{z+2}=\frac{t+2}{z+2}$; choose a rational $b$ in it and a rational
$\delta$ with $0<\delta<\min\bigl(1,\frac{(z+2)b-(t+2)}{(t+2)(z+2)b}\bigr)$. Put
$\theta=1-\delta$, $A=1-b\delta$. Then $\theta<A<1$ and
$\nu(A-\theta)=\nu(1-b)\delta>\delta=1-\theta$. By Bernoulli's inequality,
$\theta^{t+2}\ge1-(t+2)\delta$ and $A^{z+2}(1+(z+2)b\delta)\le(1-b^2\delta^2)^{z+2}\le1$;
as $(1-(t+2)\delta)(1+(z+2)b\delta)=1+\delta\bigl((z+2)b-(t+2)-(t+2)(z+2)b\delta\bigr)>1$,
\[
 A^{z+2}<\theta^{t+2}.
\]
Let $c_0=(A/\theta)^{1/s}>1$. Then $Ac_0^{t+2}<1$ (equivalent to $A^{z+2}<\theta^{t+2}$)
and $c_0^z\theta<1$ (equivalent to $A^z<\theta^t$, which follows since
$A^z=A^{z+2}/A^2<\theta^{t+2}/A^2<\theta^t$). By continuity there is a rational $c>c_0$
with $Ac^{t+2}<1$ and $c^z\theta<1$. Then $Ac^t<c^{-2}$, and $Ac^t<c^z\theta$ because
$c^s>A/\theta$. Choose a rational $B$ with $Ac^t<B<\min(c^{-2},c^z\theta)$.
\end{proof}

For $(z,t)=(1,0)$ this is \cite[Section~11.4]{L} ($A^3<\theta^2$); for $(z,t)=(2,1)$ the
key inequality is $A^4<\theta^3$, and $b\in(\frac34,1-\frac1\nu)$ needs $\nu>4$.

\begin{proof}[Proof of Theorem~\ref{I:thm:B}]
Average Proposition~\ref{I:prop:upper} over the embeddings and compare with
Lemma~\ref{I:lem:arith}, along the degrees $H$ with non-zero minors. Since
$\nu(A(1-\eta)-\bar b)-(1-\bar b)\ge g$ (the left side decreases in $\bar b<\theta$) and
$1-\bar b\le1$, the two bounds are inconsistent as soon as, with $\mathcal E=E_{\rm ar}+E_{\rm an}$,
\begin{equation}\label{I:eq:gap}
 \mathcal E<g,\qquad 1+\mathcal E<L' .
\end{equation}
We arrange $\mathcal E<\epsilon=g/2$ and $L'>1+\epsilon$. Split $\mathcal E=E_1+E_2$ with
\[
 E_1=\frac{\Lambda F_0+\log2}v+\frac{2e\log D+R}{w_0},\qquad
 E_2=\frac{\log B_R+\nu+\log2}{w_*}.
\]

\emph{Choice of $m=2n$.} Let $n$ be so large that $c^{2n}\ge2$; then
$c^{2n}/2\le Q\le c^{2n}$ and $K=Q^z$. As $n\to\infty$:
\begin{itemize}
\item $v=2S\theta^nB^{-n}\ge2Q^{z/2}\theta^nB^{-n}\ge2^{1-z/2}(c^z\theta/B)^{n}\to\infty$,
because $B<c^z\theta$; so $(\Lambda F_0+\log2)/v\to0$;
\item $(2e\log D+R)/w_0\le\bigl(2c_E\log D+2c_R+2\bigr)QB^n\le
\bigl(2c_E\log D+2c_R+2\bigr)(c^2B)^{n}\to0$, because $c^2B<1$;
\item $w_0/v=1/(2S\theta^n)$ exceeds $4(m+3)/\sigma=400(m+3)^2/(1-\theta)$ for large $n$,
and $(w_0/v)\theta^{m-1}\le\frac12$, because $c^z\theta<1$; with
$K(w_0/v)^2\theta^m\le\frac14$, $(1+3\sigma)^{k}\le\frac43$ for $k\le m+2$, and
$(1+\sigma)\theta<1$, all conditions \eqref{I:eq:volume} hold for large $n$;
\item since $v\le3Q^{z/2}\theta^nB^{-n}$ for large $n$ and $N\le c_Nc^{2tn}$,
\[
 \frac{(L')^2}{(c_2\eta^{3/2})^2}=\frac{K\theta^m}{(m+1)(m+2)v^2A^mN}
 \ \ge\ \frac{1}{9c_N(m+1)(m+2)}\Bigl(\frac{B^2}{A^2c^{2t}}\Bigr)^{n}\ \to\ \infty,
\]
because $B>Ac^t$.
\end{itemize}
Fix such an $n$ with $E_1<\epsilon/2$, $L'>1+\epsilon$ and \eqref{I:eq:volume}.

\emph{Choice of the approximations.} With $m$, $Q$ and the centers fixed, $R$ and $B_R$
are constants; choose $X$ with $(\log B_R+\nu+\log2)/X'<\epsilon/2$ for all $X'\ge X+1$,
and the $r_i$ with $w_i>X+1$ and above the thresholds. Then $\mathcal E<\epsilon$, which gives
\eqref{I:eq:gap}. This contradicts the existence of the non-zero minors $\Delta_H$ for
arbitrarily large $H$. So only finitely many good approximations exist for every
$\nu>1+(t+2)/s$, and $\mu(\xi)\le1+(t+2)/s$. The order of choices is
\[
 \nu;\ (\theta,A,B,c,\eta,\epsilon,F_0);\ n,Q;\ \text{centers},\ w_0,v,\sigma,R;\
 \text{thresholds};\ X;\ (r_i,w_i)_{i\le m};\ H\to\infty .\qedhere
\]
\end{proof}

\part{\texorpdfstring{$n$}{n} logarithms}\label{part:II}

\section{Inputs from \texorpdfstring{\cite{C}}{[C]}}\label{II:sec:inputs}

We use the homogeneous layer (files \lean{NLog*}) of the Lean library of \cite{C}, whose
centers lie on the torus, and the interpolation theorem of \cite{C}.

\paragraph{Good vectors and supplies.}
For $r\in\Q^k$ let $q(r)\ge1$ be the least common multiple of the reduced denominators of
$r_0,\dots,r_{k-1}$ (\lean{vecDen}). For $L\in\R^k$, $G\in\R$ and $\nu\in\R$, a vector
$r\in\Q^k$ is \emph{$\nu$-good} (\lean{GoodVector L G nu r}) if
\[
 |r_j|\le1\quad(j<k)\qquad\text{and}\qquad\Bigl|G-\sum_{j<k}r_jL_j\Bigr|\le q(r)^{-\nu}.
\]
$(L,G)$ is \emph{independent} (\lean{LogIndependent}) if $\sum_jy_jL_j+zG=0$ with
$y\in\Z^k$, $z\in\Z$ only for $y=0$, $z=0$. A vector $p\in\Q^k$ is \emph{off} a finite list of
pairs $(c,c_*)\in\Q^k\times\Q$ (\lean{OffHyperplanes}) if $\sum_jc_jp_j+c_*\ne0$ for every
pair in it with $c\ne0$. For a rule $\Phi$ assigning to every finite list $\mathcal L$ of
vectors a finite list $\Phi(\mathcal L)$ of such pairs, a \emph{$\Phi$-supply} of $\nu$-good
vectors (\lean{NLogSelection.Supply}) exists if for every list $\mathcal L$ and every $Q$ there
is a $\nu$-good $r$ with $q(r)\ge Q$ off $\Phi(\mathcal L)$. The rule $\Phi_k$
(\lean{NLogGeneralPosition.forbidden k}) is that of \cite[Section~2.4]{C}; for a finite list
$\mathcal H$ we write $\Phi_k^{\mathcal H}(\mathcal L)=\Phi_k(\mathcal L)\cup\mathcal H$.
\emph{General position} of $r_0,\dots,r_{m-1}$ is \cite[Definition~2.8]{C}
(\lean{NLogInterpolation.GeneralPosition}).

\begin{lemma}[Selection; \lean{exists_general_position_above_thresholds_extra}]\label{II:lem:sel}
Suppose a $\Phi_k^{\mathcal H}$-supply of $\nu$-good vectors exists. Let $m\ge0$,
$X\in\R$ and thresholds $T_0,\dots,T_{m-1}$ be given. There are $\nu$-good vectors
$r_0,\dots,r_{m-1}$ in general position and off $\mathcal H$, with $q_i=q(r_i)\ge2$,
$w_i=\lceil\log q_i\rceil>X$ and $w_i>T_i(w_0,\dots,w_{i-1})$ for all $i$.
\end{lemma}

\begin{proof}
This is \cite[Lemma~2.11]{C} (\lean{exists_above_thresholds},
\lean{generalPosition_of_admissible}) for vectors, with the list $\mathcal H$ carried along:
a vector off $\Phi_k^{\mathcal H}(\mathcal L)$ is off $\Phi_k(\mathcal L)$ and off $\mathcal H$
(\lean{offHyperplanes_append}).
\end{proof}

\paragraph{Data, characters and interpolation.}
\emph{Data} $D=(m,K,\omega,v,\theta,\sigma)$ for $k$ slopes, the margins
$\mathrm{Mar}_k$, the weights $W$ (degree) and $V$ (jets), thresholds, and weights above
thresholds are as in \cite[Section~3]{C} (\lean{NLogInterpolation.BaseData k}). For
$r\in\Q^k$ let $N(r)=\prod_j\den(r_j)$ and $p_j(r)=r_jN(r)\in\Z$ (\lean{charDen},
\lean{charExp}); the \emph{character} of $r$ is $\chi_r(Y)=Y_k^{N(r)}\prod_{j<k}Y_j^{-p_j(r)}$
(\lean{characterValue}). Centers $b_0,\dots,b_{K-1}\in(\C^\times)^{k+1}$ have \emph{distinct
characters} for $r_0,\dots,r_{m-1}$ (\lean{CharacterDistinct}) if for every $i$ the values
$\chi_{r_i}(b_c)$ are pairwise distinct.

\begin{theorem}[{\cite[Theorem~1.4]{C}}]\label{II:thm:U}
Let $k\ge1$. For all data $\mathcal D$ for $k$ slopes there are thresholds $T(\mathcal D)$ such that for all
rational weights $w$ above $T(\mathcal D)$ there is $R_0\ge1$ with the following property. For all
slopes $r_0,\dots,r_{m-1}\in\Q^k$ in general position and all centers
$b_0,\dots,b_{K-1}\in(\C^\times)^{k+1}$ with distinct characters there is $H_0$ such that the
packet map of the polynomials of $W$-degree at most $H$ to the jets of $V$-weight less than $H$
at the centers, with the chart
\[
 Y_j\mapsto b_{c,j}(1+t_j),\qquad X_i\mapsto u_i+\log(1+t_k)-\sum_{j<k}r_{i,j}\log(1+t_j),
\]
is surjective for every $H\ge H_0$ divisible by $R_0$.
\end{theorem}

In Lean this is \lean{NLogInterpolation.UniformStatement k}, obtained from
\lean{CLogUniformInterpolation.cLogUniformInterpolation} by
\lean{CLogInterpolation.uniformStatement_of_C}. It is the special case of \cite[Theorem~1.4]{C} for centers on the torus, without additive
coordinates, and it holds for arbitrary such centers; \cite{C} applies it to the centers
$\beta^J$.

\paragraph{Other inputs.}
From \cite{C}: the truncation lemma \cite[Lemma~6.1]{C}, the counts of rows and the mean
$\overline b_H\in[0,\theta]$ \cite[(6.3)]{C}, the collision count
$\sum|p|\ge\frac14N^{\frac{k+2}{k+1}}-N$ for $N$ distinct $p\in\Z_{\ge0}^{k+1}$
\cite[Lemma~7.1]{C}, the growth and row-test bounds of \cite[Section~7.3]{C}, the occupancy
dichotomy of \cite[Section~7.5]{C}, and the margins lemma \cite[Lemma~8.3]{C}
(\lean{NLogInterpolation.margins_of_pos_bounds}). From \cite{L}: $\log\operatorname{lcm}(1,\dots,n)\le\Lambda_0n$ with
$\Lambda_0=4+\log4$, the norm comparison over a number field, and the comparison of
determinant bounds along cofinally many degrees (\lean{CofinalComparison}). From Part~I:
boxes of products and torsion coordinates (\lean{TwoGridPairsBasic}; the indexing
\lean{TwoGrid.boxIdx}, \lean{TwoGrid.boxCenter} of \lean{TwoGridPairs}), the description of the
embeddings of a Galois field (\lean{TwoGrid.exists_aut_of_embedding}), two elementary
inequalities (\lean{TwoGrid.bernoulli_gap_pq}, \lean{TwoGrid.lt_one_of_pow_lt_one}) and the
structure \lean{TwoGrid.GridParams} of \lean{TwoGridData}. No theorem of Part~I on the
irrationality exponent is used.

\section{Reduction and the induction}\label{II:sec:red}

Let $\Meas(n)$ be the statement of Theorem~\ref{II:thm:main} in the form of Lean
(\lean{MeasureAlg n}): for every Galois number field $E$, embedding $\iota:E\to\C$ and
$\alpha\in E^n$ with $\iota(\alpha_i)=a_i>0$ without a power relation, and every family
$s_1,\dots,s_z\in\Gal(E/\Q)$ which is \emph{independent} (\lean{TupleIndep}:
$\prod_l s_l(\alpha_i)^{e_l}=1$ for all $i$ only for $e=0$) and of the largest size among
independent families, the bound \eqref{II:eq:main} holds with $z$ in place of $\rho(a)$.
A largest family exists (\lean{exists_largest}), and $z=\rho(a)$ for $E$ the Galois closure.
$\Meas(1)$ holds trivially: $|h\log a|\ge|\log a|$ (\lean{measureAlg_one}).

\begin{proposition}[\lean{measureAlg_succ_of_bounded}]\label{II:prop:bounded}
Let $k\ge1$. Suppose that for every Galois number field $E$, $\iota$, $\alpha\in E^k$,
$\gamma\in E$ with $\iota(\alpha_j)=a_j>0$, $\iota(\gamma)=g>0$ and $(L,G)$ independent, for
every largest independent family of conjugate tuples of $(\alpha,\gamma)$, of size $z$, and
every $\nu>k+1+z$, the $\nu$-good vectors for $(L,G)$ have bounded denominators. Then
$\Meas(k+1)$ holds.
\end{proposition}

\begin{proof}
Let $\alpha\in E^{k+1}$, $a=\iota(\alpha)$, $s$ and $z$ be as in $\Meas(k+1)$, and
$\varepsilon>0$. For an index $i_0$ consider the reordered tuple
$((\alpha_i)_{i\ne i_0},\alpha_{i_0})$. Its conjugate tuples are those of $\alpha$ with the
coordinates permuted, so $s$ is again independent and largest (\lean{tupleIndep_of_cover}),
and the logarithms are independent since $a$ has no power relation. Apply the hypothesis with
$\nu=k+z+1+\varepsilon$. As in \cite[Lemma~2.2]{C} (\lean{ordered_bound_shift}), with
$\kappa=k+z-1$ in place of $k$: if $\delta>0$ bounds $|G-\sum_jr_jL_j|$ from below on the
finite set of $r$ with $|r_j|\le1$ and $q(r)$ below the bound, then $C=\min(1,\delta)$
satisfies
\[
 C|w|^{-(k+z)-\varepsilon}\le\Bigl|\sum_jx_jL_j+wG\Bigr|\qquad(w\ne0,\ |x_j|\le|w|),
\]
since otherwise $r=-x/w$ would be a good vector with $q(r)\le|w|$ and
$|G-\sum_jr_jL_j|<\delta$. Taking the smallest of the constants over $i_0$, and $i_0$ with
$|h_{i_0}|$ largest, gives \eqref{II:eq:main} with the exponent $-(k+z)-\varepsilon$, that is
$-((k+1)-1+z)-\varepsilon$.
\end{proof}

\begin{lemma}[\lean{tupleIndep_of_reduced}, \lean{tupleIndep_of_monomial}]\label{II:lem:mono}
Let $\tilde\alpha_j=\prod_i\alpha_i^{M_{ji}}$ be monomials in the coordinates of $\alpha$. A
family of automorphisms independent for $\tilde\alpha$ is independent for $\alpha$. Hence the
rank of the conjugate tuples of $\tilde\alpha$ is at most that of $\alpha$
(\lean{exists_largest_monomial}).
\end{lemma}

\begin{proof}
If $\prod_ls_l(\alpha_i)^{e_l}=1$ for all $i$, then
$\prod_ls_l(\tilde\alpha_j)^{e_l}=\prod_i\bigl(\prod_ls_l(\alpha_i)^{e_l}\bigr)^{M_{ji}}=1$ for
all $j$, so $e=0$.
\end{proof}

\begin{lemma}[Hyperplanes; \lean{hyperplaneBounded_alg}]\label{II:lem:plane}
Assume $\Meas(k)$. Let $E,\iota,\alpha,\gamma,a,g$ be as in Proposition~\ref{II:prop:bounded},
with $(L,G)$ independent, let $z$ be the size of a largest independent family for
$(\alpha,\gamma)$, and $\nu>k+1+z$. For $c\in\Q^k\setminus\{0\}$ and $c_*\in\Q$ the
$\nu$-good vectors on the hyperplane $\sum_jc_jr_j+c_*=0$ have bounded denominators.
\end{lemma}

\begin{proof}
We may assume $c\in\Z^k$, $c_*\in\Z$; choose $i$ with $c_i\ne0$. As in
\cite[Lemma~2.5]{C}, put
\[
 \Lambda_j=c_iL_j-c_jL_i\ (j\ne i),\qquad\Lambda_i=c_iG+c_*L_i ;
\]
these are the logarithms of $\tilde a_j=a_j^{c_i}a_i^{-c_j}$ and $\tilde a_i=g^{c_i}a_i^{c_*}$
(\lean{reducedBaseR}), which are the images under $\iota$ of the monomials
$\tilde\alpha_j=\alpha_j^{c_i}\alpha_i^{-c_j}$, $\tilde\alpha_i=\gamma^{c_i}\alpha_i^{c_*}$ in
$E$ (\lean{reducedBaseE}). They have no power relation, because in a relation
$\sum_jy_j\Lambda_j=0$ the coefficient of $G$ is $y_ic_i$ and that of $L_j$, $j\ne i$, is
$y_jc_i$. By Lemma~\ref{II:lem:mono} the rank $\tilde z$ of their conjugate tuples is at most $z$.
$\Meas(k)$ with $\varepsilon=\frac12$ gives $C>0$ with
$|\sum_jx_j\Lambda_j|\ge C\,H(x)^{-(k-1+\tilde z)-1/2}$ for $x\in\Z^k\setminus0$. For a good
$r$ on the hyperplane with $q=q(r)$ let $x_j=-qr_j$ $(j\ne i)$, $x_i=q$; then
$\sum_jx_j\Lambda_j=c_iq\,(G-\sum_jr_jL_j)$ and $H(x)\le q$, so
\[
 C\,q^{-(k-1+\tilde z)-1/2}\le|c_i|\,q^{1-\nu}\le|c_i|\,q^{-(k+\tilde z)}
\]
because $\nu\ge k+1+\tilde z$; hence $q\le(|c_i|/C)^2$ (\lean{hyperplane_bound_core_shift}).
\end{proof}

\begin{proposition}[\lean{measureAlg_succ_of_no_supply}]\label{II:prop:nosupply}
Let $k\ge1$ and assume $\Meas(k)$. Suppose that, in the setting of
Proposition~\ref{II:prop:bounded}, for every $\nu>k+1+z$ it is impossible that a
$\Phi_k^{\mathcal H}$-supply of $\nu$-good vectors exists for every finite list $\mathcal H$.
Then $\Meas(k+1)$ holds.
\end{proposition}

\begin{proof}
By Proposition~\ref{II:prop:bounded} it suffices to bound the denominators of the good vectors.
If they are unbounded, then for every list $\mathcal L$ Lemma~\ref{II:lem:plane}, applied to the
finitely many pairs of $\Phi_k^{\mathcal H}(\mathcal L)$ with $c\ne0$, shows that good vectors
of large denominator are off $\Phi_k^{\mathcal H}(\mathcal L)$
(\lean{exists_bound_offHyperplanes_of}); so a $\Phi_k^{\mathcal H}$-supply exists for every
$\mathcal H$ (\lean{supply_of_unbounded_of}).
\end{proof}

\begin{theorem}[No supply; \lean{no_supply_alg}]\label{II:thm:nosupply}
Let $k\ge1$, let $E,\iota,\alpha,\gamma,a,g$ be as in Proposition~\ref{II:prop:bounded} with
$(L,G)$ independent, let $s_1,\dots,s_z$ be a largest independent family of conjugate tuples
of $(\alpha,\gamma)$, and $\nu>k+1+z$. Then there is a finite list $\mathcal H$ for which no
$\Phi_k^{\mathcal H}$-supply of $\nu$-good vectors for $(L,G)$ exists.
\end{theorem}

The list will be the list $\mathcal H_Q$ of Proposition~\ref{II:prop:supply}(T5$'$) for a value
of $Q$ determined by the parameters. Theorem~\ref{II:thm:nosupply} is proved in
Section~\ref{II:sec:end}; it does not use $\Meas(k)$.

\begin{proof}[Proof of Theorem~\ref{II:thm:main}]
$\Meas(1)$ holds, and $\Meas(k)$ implies $\Meas(k+1)$ for $k\ge1$ by
Proposition~\ref{II:prop:nosupply} and Theorem~\ref{II:thm:nosupply} (\lean{measureAlg_succ}). So
$\Meas(n)$ holds for all $n\ge1$ (\lean{measureAlg_all}); for $E$ the Galois closure of
$\Q(a)$ and a largest family, $z=\rho(a)$.
\end{proof}

\section{Tuple families}\label{II:sec:supply}

\begin{definition}[\lean{TupleFamily E a t}]\label{II:def:family}
Let $E$ be a number field and $a=(a_0,\dots,a_k)$ positive reals. A \emph{tuple family of
size $z=t+1$} for $a$ consists of $\beta_0,\dots,\beta_t\in(E^\times)^{k+1}$ such that
\begin{enumerate}
\item[(i)] the $\beta_l$ are independent: $\prod_l\beta_{l,\tau}^{e_l}=1$ for all $\tau$ only
for $e=0$;
\item[(ii)] $D\beta_{l,\tau}$ is an algebraic integer for all $l,\tau$, for some integer
$D\ge1$;
\item[(iii)] for every embedding $\sigma:E\to\C$ there are $A\in\Z$ and
$e\in\Z^{t+1}\setminus0$ with $\prod_l\sigma(\beta_{l,\tau})^{e_l}=a_\tau^{A}$ for all
$\tau\le k$.
\end{enumerate}
\end{definition}

\begin{lemma}[Relation; \lean{exists_relation_of_maximal_tuple}]\label{II:lem:relation}
Let $s_1,\dots,s_z$ be a largest independent family of conjugate tuples of $\alpha$. For
every $g\in\Gal(E/\Q)$ there are $A\in\Z\setminus0$ and $e\in\Z^z$ with
$\prod_ls_l(\alpha_i)^{e_l}=g(\alpha_i)^A$ for all $i$.
\end{lemma}

\begin{proof}
The family $g,s_1,\dots,s_z$ is not independent. In a non-trivial relation the exponent of
$g$ is not zero, since $s_1,\dots,s_z$ are independent; its negative is $A$.
\end{proof}

\begin{proposition}[\lean{exists_tupleFamily}]\label{II:prop:family}
Let $E$ be Galois, $\iota:E\to\C$, $\alpha\in E^{k+1}$ with $\iota(\alpha)=a\in\R_{>0}^{k+1}$
without a power relation, and $s_1,\dots,s_z$ a largest independent family. Then $z\ge1$ and
the tuples $\beta_l=s_l(\alpha)$ form a tuple family of size $z$ for $a$.
\end{proposition}

\begin{proof}
The family consisting of the identity alone is independent: $\alpha_0^e=1$ gives
$a_0^e=1$ and $e=0$. So $z\ge1$. (i) is the independence of $s$. (ii): if $D$ is a common
denominator of the $\alpha_i$, then $Ds_l(\alpha_i)=s_l(D\alpha_i)$ is integral. (iii): every
embedding is $\sigma=\iota\circ g$ with $g\in\Gal(E/\Q)$ (\lean{exists_aut_of_embedding}).
Lemma~\ref{II:lem:relation} for $g^{-1}$ gives $\prod_ls_l(\alpha)^{e_l}=g^{-1}(\alpha)^A$, and
applying $\iota\circ g$, $\prod_l\sigma(\beta_l)^{e_l}=a^A$. If $e=0$, then
$g^{-1}(\alpha_0)^A=1$, hence $a_0^A=1$ and $A=0$, which is excluded.
\end{proof}

For an integer $Q\ge1$ the \emph{centers} of a tuple family are (\lean{TupleFamily.center})
\[
 b_c=\prod_l\beta_l^{c_l}\in(E^\times)^{k+1},\qquad c\in\{0,\dots,Q-1\}^{z},
\]
$K=Q^z$ of them.

\begin{proposition}[The supply]\label{II:prop:supply}
Let $\beta_0,\dots,\beta_t$ be a tuple family of size $z=t+1$ for $a$, and
$L=(\log a_0,\dots,\log a_k)$. There are integers $N_b\ge0$ and a real $c_R\ge0$ such that,
with $c_N=(N_b+1)(2N_b+3)^t$:
\begin{enumerate}
\item[(T2)] $D^{zQ}b_{c,\tau}$ is an algebraic integer for all $c,\tau$
(\lean{center_integral});
\item[(T3, T4$'$)] for every $Q\ge1$ and every embedding $\sigma$ there are vectors
$Y_c\in\C^{k+1}$, a map $c\mapsto\mu(c)$ from the box to $\{0,\dots,N-1\}$ with $N=c_NQ^t$, vectors
$Y^0_0,\dots,Y^0_{N-1}\in\C^{k+1}$ and integers $s_c$ such that
\[
 e^{Y_{c,\tau}}=\sigma(b_{c,\tau}),\qquad Y_c=Y^0_{\mu(c)}+s_cL,\qquad
 |Y_{c,\tau}|\le c_RQ,\quad|Y^0_{\mu,\tau}|\le c_R(N_b+1)Q,\quad|s_c|\le c_RQ
\]
(\lean{exists_logs});
\item[(T5$'$)] for every $Q$ there is a finite list $\mathcal H_Q$ of pairs
$(c,c_*)\in\Q^k\times\Q$ such that for every embedding $\sigma$ and every $r\in\Q^k$ off
$\mathcal H_Q$ the values $\chi_r(\sigma(b_c))$ are pairwise distinct
(\lean{TupleFamily.exists_hyperplanes}).
\end{enumerate}
\end{proposition}

\begin{proof}
(T2): $D^{zQ}b_{c,\tau}=D^{zQ-\sum_lc_l}\prod_l(D\beta_{l,\tau})^{c_l}$.

(T3, T4$'$). Fix $\sigma$ and the relation (iii), with signs changed so that $e_{l_0}>0$ for
some $l_0$ (\lean{exists_rel_log}). There are integers $\kappa_\tau$ with
\[
 \sum_le_l\log\sigma(\beta_{l,\tau})=AL_\tau+2\pi i\,\kappa_\tau\qquad(\tau\le k).
\]
For $c$ in the box let $j=\lfloor c_{l_0}/e_{l_0}\rfloor$, $c^*=c-je$ and
\[
 Y_{c,\tau}=\sum_lc_l\log\sigma(\beta_{l,\tau})-2\pi i\,j\kappa_\tau
 =\sum_lc^*_l\log\sigma(\beta_{l,\tau})+jA\,L_\tau
\]
(\lean{box_log_eq}). So $e^{Y_c}=\sigma(b_c)$, $Y_c=Y^0(c^*)+s_cL$ with $s_c=jA$ and
$Y^0(c^*)_\tau=\sum_lc^*_l\log\sigma(\beta_{l,\tau})$, and the label of $c$ is $c^*$. Let
$N_b$ bound $|e_l|$ for all $\sigma$ and $l$. Then $0\le c^*_{l_0}<e_{l_0}$ and
$|c^*_l|\le(N_b+1)Q$, so $c^*_{l_0}$ takes at most $N_b+1$ values and every other coordinate
at most $(2N_b+3)Q$, and $c^*$ takes at most $c_NQ^t$ values (\lean{card_image_reduced_le});
enumerate them (\lean{exists_labels_fn}). The bounds hold if $c_R$ bounds
$\sum_{\tau}\bigl(\sum_l|\log\sigma(\beta_{l,\tau})|+2\pi|\kappa_\tau|\bigr)+|A|$ for every
embedding.

(T5$'$) is Lemma~\ref{II:lem:T5} below for the units $\beta_{l,\tau}$: since $\chi_r$ is a
homomorphism and $\sigma$ is injective, $\chi_r(\sigma(b_c))=\chi_r(\sigma(b_{c'}))$ means
$\chi_r(u)=1$ for $u=\prod_l\beta_l^{c_l-c'_l}$, and $|c_l-c'_l|\le Q$
(\lean{characterDistinct_box}).
\end{proof}

\begin{lemma}[\lean{exists_hyperplanes}]\label{II:lem:T5}
Let $M$ be an abelian group, written additively, and $P_0,\dots,P_t\in M^{k+1}$ independent:
$\sum_le_lP_{l,\tau}=0$ for all $\tau$ only for $e=0$. For every $Q$ there is a finite list
$\mathcal H_Q$ of pairs $(c,c_*)\in\Q^k\times\Q$ such that for $r\in\Q^k$ off $\mathcal H_Q$ and
$\Delta\in\Z^{t+1}$ with $|\Delta_l|\le Q$,
\[
 N(r)\,x_k-\sum_{j<k}p_j(r)\,x_j=0,\qquad x=\sum_l\Delta_lP_l ,
\]
holds only for $\Delta=0$.
\end{lemma}

\begin{proof}
Let $\psi_1,\dots,\psi_g:\langle P_{l,\tau}\rangle\to\Z$ be homomorphisms whose common kernel is
the torsion subgroup (\lean{exists_torsion_coordinates}). For $\Delta$ in the box and $i\le g$
let $\mathcal H_Q$ contain the pair
\[
 c_j=-\psi_i(x_j)\ (j<k),\qquad c_*=\psi_i(x_k),\qquad x=\sum_l\Delta_lP_l .
\]
If the displayed relation holds, applying $\psi_i$ and dividing by $N(r)$ gives
$\psi_i(x_k)-\sum_jr_j\psi_i(x_j)=0$, that is, $r$ lies on the hyperplane of $(\Delta,i)$. As
$r$ is off $\mathcal H_Q$, the linear part vanishes: $\psi_i(x_j)=0$ for $j<k$, and then
$\psi_i(x_k)=0$ by the relation. So all $x_\tau$ are torsion, $\sum_lN\Delta_lP_{l,\tau}=0$ for
some $N\ge1$, and $\Delta=0$ by independence.
\end{proof}

For $k=1$ the hyperplanes are points of bounded denominator; this is the condition (T5) of
Definition~\ref{I:def:supply}.

\section{The determinant and its arithmetic lower bound}\label{II:sec:ar}

\paragraph{The setting.}
Throughout Sections~\ref{II:sec:ar}--\ref{II:sec:an} the following are fixed (structure
\lean{NAlgData k nu}): real numbers $L=(L_0,\dots,L_k)$ with $L_0\ne0$; parameters
$0<\theta<A<1$ ($\theta,A$ rational) and $\eta\in(0,1)$ with
\begin{equation}\label{II:eq:gap}
 \mathrm{gap}:=\nu\bigl(A(1-\eta)-\theta\bigr)-(1-\theta)>0
\end{equation}
(\lean{TwoGrid.GridParams nu}); data $D=(m,K,\omega,v,\theta,\sigma)$ for $k$ slopes with
the same $\theta$; slopes $r_0,\dots,r_{m-1}\in\Q^k$ in general position with $|r_{i,j}|\le1$,
$q_i=q(r_i)\ge2$ and
\[
 \Bigl|L_k-\sum_{j<k}r_{i,j}L_j\Bigr|\le q_i^{-\nu};
\]
the weights $w_i=\lceil\log q_i\rceil$, a number $w_*>0$ with $w_*\le w_i$, a radius $R\ge1$,
and a number $N\ge1$ of labels.

\paragraph{The matrix.}
Let $b_0,\dots,b_{K-1}\in(E^\times)^{k+1}$ be centers in a number field $E$. With
$F_0=4/\theta$, $T_i=\lceil F_0w_i/v\rceil$, $G_i(x)=\sum_{1\le l<T_i}(-1)^{l+1}x^l/l$, and
the rows $(c,\lambda,\beta)$ and columns $Y^hX^\gamma$ of \cite[Section~6.2]{C}, the matrix
$A_H$ (\lean{algMatrix}) has the entries
\begin{equation}\label{II:eq:matrix}
 [t^\lambda u^\beta]\ \prod_{j\le k}\bigl(b_{c,j}(1+t_j)\bigr)^{h_j}
 \prod_{i<m}\Bigl(u_i+G_i(t_k)-\sum_{j<k}r_{i,j}G_i(t_j)\Bigr)^{\gamma_i} .
\end{equation}
This is the matrix of \cite{C} with the center $\beta^J$ replaced by $b_c$: in Lean it is the
entry of \cite{C} with $J=1$ and the base $b_c$, so that every lemma on that entry applies. Let
$M=M_H$ be the number of rows; a selection of $M$ columns defines a minor $\Delta_H$
(\lean{algSelectedMinor}). For an embedding $\tau$, $\tau(\Delta_H)$ is the minor built with
the centers $\tau(b_c)$ (\lean{map_algSelectedMinor_det}).

\begin{lemma}[Non-zero minors; \lean{cofinal_field_alg_minors}]\label{II:lem:minors}
Let $T(\mathcal D)$ be the thresholds of Theorem~\ref{II:thm:U} and suppose the weights are above $T(\mathcal D)$.
If for one embedding $\sigma_0$ the centers $\sigma_0(b_c)$ have distinct characters for
$r_0,\dots,r_{m-1}$, then for every $H_1$ there are an integer $H\ge H_1$, $H>0$, and a
selection with $\Delta_H\ne0$; then $\tau(\Delta_H)\ne0$ for every embedding $\tau$.
\end{lemma}

\begin{proof}
Theorem~\ref{II:thm:U} applies to the centers $\sigma_0(b_c)\in(\C^\times)^{k+1}$; the
full-logarithm matrix has rank $M$ for large $H$ divisible by $R_0$, and by the truncation
lemma \cite[Lemma~6.1]{C} the same minor of \eqref{II:eq:matrix} at $\sigma_0(b)$ is
$\sigma_0(\Delta_H)\ne0$ (\lean{cofinal_alg_minors}). Embeddings are injective.
\end{proof}

\begin{proposition}[Arithmetic lower bound; \lean{algSelectedMinor_arithmetic_lower_bound}]\label{II:prop:lower}
Let $d\ge1$ be an integer with $d\,b_{c,t}$ integral for all $c,t$, and put
\[
 E_{\rm ar}(d)=\frac{\Lambda_0F_0}v+\frac{(k+1)\log d}{\omega}
\]
(\lean{arithmeticError}). For every $H>0$ and every selection with $\Delta_H\ne0$,
\[
 \frac1{[E:\Q]}\sum_{\tau:E\to\C}\frac{\log|\tau(\Delta_H)|}{M_HH}\ \ge\
 -(1-\overline b_H)-E_{\rm ar}(d).
\]
\end{proposition}

\begin{proof}
As in \cite[Proposition~6.3]{C}, with the denominator of the bases replaced by
$d^{(k+1)\lfloor H/\omega\rfloor}$ (\lean{algDenominator}): since $h_j\le H/\omega$ for every
column, this integer times $\prod_jb_{c,j}^{h_j}$ is an algebraic integer. Multiplying the
columns, rows and entries as there makes all entries algebraic integers
(\lean{algSelectedMinor_entries_integral}), so the cleared minor is a non-zero algebraic
integer $x$ and $\sum_\tau\log|\tau(x)|=\log|N_{E/\Q}(x)|\ge0$
(\lean{algSelectedMinor_clearing_bound}). The costs are those of \cite{C}
(\lean{alg_combined_log_cost_le}) and $\log d^{(k+1)\lfloor H/\omega\rfloor}\le(k+1)H\log d/\omega$
(\lean{log_algDenominator_le}).
\end{proof}

\begin{corollary}[\lean{algSelectedMinor_uniform_upper_bound_obstruction}]\label{II:cor:obs}
If moreover $\log|\tau(\Delta_H)|\le M_HH\,U$ for every embedding $\tau$, then
$-(1-\overline b_H)-E_{\rm ar}(d)\le U$.
\end{corollary}

\section{The analytic upper bound with labels}\label{II:sec:an}

\paragraph{Exact forms.}
As in \cite[Section~7.1]{C} with constants $0$: let
$\eta_i=(L_k-\sum_jr_{i,j}L_j)/L_0$ and $\ell^*_i(z)=z_k-\sum_jr_{i,j}z_j-\eta_iz_0$
(\lean{eta}, \lean{exactCoeffR}). Then
\begin{equation}\label{II:eq:exact}
 \ell^*_i(L)=0,\qquad|\eta_i|\le\frac{q_i^{-\nu}}{|L_0|},\qquad
 \sum_t|c^*_{i,t}|\le M_L:=k+1+\frac1{|L_0|}
\end{equation}
(\lean{exact_relation}, \lean{abs_eta_le}, \lean{slope_norm_bound}). For a column
$P=Y^hX^\gamma$ and $a\le\gamma$ put
\[
 f_{a,P}(z)=e^{h\cdot z}\prod_i\binom{\gamma_i}{a_i}\ell^*_i(z)^{\gamma_i-a_i}
 \qquad(z\in\C^{k+1}).
\]

\paragraph{Labelled centers.}
At one embedding the centers are given with their logarithms (structure \lean{Centers}):
points $b_c\in(\C^\times)^{k+1}$, a label map $\mu:\{0,\dots,K-1\}\to\{0,\dots,N-1\}$,
representatives $Y^0_\mu\in\C^{k+1}$ and integers $s_c$ with
\begin{equation}\label{II:eq:centers}
 \exp\bigl(Y^0_{\mu(c),t}+s_cL_t\bigr)=b_{c,t},\qquad
 |s_cL_t|+\tfrac34\le\tfrac R2,\qquad|Y^0_{\mu,t}|\le R .
\end{equation}

\begin{lemma}[\lean{mainEntry_eq_columnRow}]\label{II:lem:shift}
For every center $c$, column $P$ and $a\le\gamma$,
\[
 [t^\lambda]\ \prod_j\bigl(b_{c,j}(1+t_j)\bigr)^{h_j}\prod_i\binom{\gamma_i}{a_i}
 \ell^*_i\bigl(\log(1+t)\bigr)^{\gamma_i-a_i}
 \ =\ e^{h\cdot Y^0_{\mu(c)}}\ [t^\lambda]\,f_{a,P}\bigl(s_cL+\log(1+t)\bigr).
\]
\end{lemma}

\begin{proof}
$\prod_jb_{c,j}^{h_j}=e^{h\cdot Y^0_{\mu(c)}}e^{h\cdot s_cL}$, and
$\ell^*_i(s_cL+\log(1+t))=\ell^*_i(\log(1+t))$ because $\ell^*_i(L)=0$
(\lean{periodSeries_center_collapse}).
\end{proof}

So the row expansion of \cite[Section~7.2]{C} holds with $z_c(t)=s_cL+\log(1+t)$ in place of
$z_J(t)=JL+\log(1+t)$ and the factor $e^{h\cdot Y^0_{\mu(c)}}$ in the column function
(\lean{det_slopeEntry_eq_sum_bases}, \lean{det_algFullLog_selected_eq_algSelectedMinor}):
\begin{equation}\label{II:eq:assign}
 \Delta_H=\sum_{\rm assignments}\ \prod_{\rm rows}s(a,\beta,\ell)\cdot
 \det\Bigl(e^{h\cdot Y^0_{\mu(c)}}\,[t^{\lambda-\ell e_0}]\,f_{a,P}\bigl(z_c(t)\bigr)\Bigr),
\end{equation}
with the scalars $s(a,\beta,\ell)$ of \cite[(7.5)]{C}, which do not depend on the centers and satisfy \cite[(7.7)]{C}.
The center $s_cL$ lies in the polydisc where the row tests of \cite[(7.9)]{C} hold, by
\eqref{II:eq:centers}.

\begin{lemma}[Labelled collisions; \lean{weighted_label_determinant_exp_bound}]\label{II:lem:det}
Fix an assignment, call $(\mu(c),a)$ the \emph{group} of a row $(c,\lambda,\beta)$ with
choice $(a,\ell)$, and let $n_\kappa$ be the number of rows in the group $\kappa$. Then the
determinant in \eqref{II:eq:assign} is at most
\[
 \exp\Bigl(-c_0\sum_\kappa n_\kappa^{\frac{k+2}{k+1}}
 +MH\Bigl(\frac{2R}\omega+\frac{\log(2RM_L)}{w_*}+\frac{\log2}v+\mathrm{rem}_H\Bigr)\Bigr)
 \prod_{\rm rows}2^{-\ell}
\]
in absolute value, with $c_0=\frac{\log2}8$ and $\mathrm{rem}_H$ of
\cite[Lemma~7.2]{C}.
\end{lemma}

\begin{proof}
Expand $f_{a,P}(z)=\sum_pc_p(a,P)\prod_t(z_t/R)^{p_t}$. The coefficient row of a row of the
group $(\mu,a)$ with index $p$ is $\bigl(e^{h\cdot Y^0_\mu}c_p(a,P)\bigr)_P$; it depends only on
$(\mu,a,p)$, so terms in which two rows of one group have the same $p$ vanish, and the proof of
\cite[Lemma~7.2]{C} applies with the groups $(\mu,a)$ in place of $a$. Rows of one group
may belong to different centers $c$, with different points $s_cL$; this does not matter,
because the center enters only the row test $[t^{\lambda-\ell e_0}]\prod_t(z_{c,t}(t)/R)^{p_t}$
and not the coefficient row. The only new factor is
$|e^{h\cdot Y^0_\mu}|\le e^{R|h|}\le e^{RH/\omega}$ per row, which together with the growth bound
$\exp(H(R/\omega+\log(2RM_L)/w_*))$ gives $2R/\omega$.
\end{proof}

Let $D_A(H)$ be the number of $a\in\Z_{\ge0}^m$ with $\sum_iw_ia_i\le AH$. The groups $(\mu,a)$
with $\sum_iw_ia_i\le AH$ are at most $N\,D_A(H)$ (\lean{low_label_capacity}), so H\"older's
inequality gives $\sum_\kappa n_\kappa^{(k+2)/(k+1)}\ge(\eta M)^{(k+2)/(k+1)}(ND_A(H))^{-1/(k+1)}$
when at least $\eta M$ rows lie in such groups. Put
\[
 \Gamma_H=\frac{c_0\,\eta^{\frac{k+2}{k+1}}}{H}
 \Bigl(\frac{M_H}{N\,D_A(H)}\Bigr)^{\frac1{k+1}}\ \longrightarrow\
 \Gamma=c_0\,\eta^{\frac{k+2}{k+1}}\Bigl(\frac{K\theta^m}{(m+1)\cdots(m+k+1)\,v^{k+1}A^m\,N}\Bigr)^{\frac1{k+1}}
\]
(\lean{collisionRate}, \lean{collisionLimit}), and
\begin{equation}\label{II:eq:Ean}
 E_{\rm an}=\frac{2R}\omega+\frac{\log2}v+\frac{\log(2RM_L)+\nu+\log2+\log c_L}{w_*},
 \qquad c_L=\max\Bigl(1,\frac1{|L_0|}\Bigr)
\end{equation}
(\lean{analyticError}).

\begin{proposition}[Analytic upper bound; \lean{actualMinor_norm_upper_bound}]\label{II:prop:upper}
Let $\nu\ge0$ and let centers be given as in \eqref{II:eq:centers}. For every integer $H>0$ and
every selection,
\[
 \frac{\log|\Delta_H|}{M_HH}\le E_{\rm an}+\mathrm{rem}'_H
 +\max\bigl\{-\Gamma_H,\ -\nu\bigl(A(1-\eta)-\overline b_H\bigr)\bigr\},
\]
where $\Delta_H$ is the minor of \eqref{II:eq:matrix} at these centers and
$\mathrm{rem}'_H\to0$ is the remainder of \cite[Proposition~7.3]{C} (\lean{analyticRemainder}).
\end{proposition}

\begin{proof}
As \cite[Proposition~7.3]{C}: the occupancy dichotomy with the groups $(\mu,a)$
(\lean{collision_and_scalar_exponent_le_data}), Lemma~\ref{II:lem:det}, the scalar bound
$|s(a,\beta,\ell)|\le2^\ell\exp(-\nu\sum_iw_i(a_i-\beta_i)+(\nu+\log2+\log c_L)H/w_*)$, and the
count of assignments (\lean{assignmentTerm_norm_le}).
\end{proof}

The bound holds at every embedding at which the centers have labelled logarithms; for the
centers of a tuple family this is every embedding, by Proposition~\ref{II:prop:supply}(T3,~T4$'$).

\section{Parameters and the end of the proof}\label{II:sec:end}

\begin{lemma}[\lean{exists_nalg_parameters}]\label{II:lem:params}
Let $p\ge1$, $t\ge0$, $s\ge1$ be integers, $z=t+s$, and $\nu s>z+p$. There are rationals
$\theta,A,B,c$ with
\begin{equation}\label{II:eq:P}
 0<\theta<A<1,\quad\nu(A-\theta)>1-\theta,\quad c>1,\quad Ac^t<B,\quad c^pB<1,\quad
 B<c^z\theta<1 .
\end{equation}
\end{lemma}

\begin{proof}
This is Lemma~\ref{I:lem:params} with the height power $2$ replaced by $p$. The interval
$\frac{t+p}{z+p}<b<1-\frac1\nu$ is non-empty since $1-\frac1\nu>1-\frac s{z+p}$; choose a
rational $b$ in it and a rational $\delta$ with
$0<\delta<\min\bigl(1,\frac{(z+p)b-(t+p)}{(t+p)(z+p)b}\bigr)$, and put $\theta=1-\delta$,
$A=1-b\delta$. Then $\theta<A<1$ and $\nu(A-\theta)=\nu(1-b)\delta>\delta=1-\theta$. By
Bernoulli's inequality $\theta^{t+p}\ge1-(t+p)\delta$ and
$A^{z+p}(1+(z+p)b\delta)\le(1-b^2\delta^2)^{z+p}\le1$, and
$(1-(t+p)\delta)(1+(z+p)b\delta)>1$ by the choice of $\delta$; so $A^{z+p}<\theta^{t+p}$
(\lean{bernoulli_gap_pq}). Let $c_0=(A/\theta)^{1/s}>1$. Then $Ac_0^{t+p}<1$ (its $s$-th power
times $\theta^{t+p}$ is $A^{z+p}$) and $c_0^z\theta<1$ (its $s$-th power times $\theta^t$ is
$A^z$, and $A^z<\theta^t$ follows from $A^{z+p}<\theta^{t+p}$ and $\theta^p<A^p$). Choose a
rational $c>c_0$ with $Ac^{t+p}<1$ and $c^z\theta<1$; then $Ac^t<c^{-p}$, and $Ac^t<c^z\theta$
because $c^s>A/\theta$. Choose a rational $B$ with $Ac^t<B<\min(c^{-p},c^z\theta)$.
\end{proof}

We use $p=k+1$, $s=1$, $t=z-1$: the condition is $\nu>k+1+z$
(\lean{exists_nalg_parameters_rank}).

\paragraph{The dimension.}
For parameters \eqref{II:eq:P} with $p=k+1$ and an integer $T\ge1$ put (\lean{algQ},
\lean{algK}, \lean{algW}, \lean{algV})
\begin{equation}\label{II:eq:values}
 m=(k+1)T,\quad Q=\lfloor c^{(k+1)T}\rfloor,\quad K=Q^z,\quad\omega=B^{-T},\quad
 v=\frac{2(c^z\theta)^T}{B^T},\quad\sigma=\frac1{100(m+k+4)} .
\end{equation}
Then $\omega/v=\frac12(c^z\theta)^{-T}$ and $(\omega/v)\theta^T=\frac12c^{-zT}$, so
\begin{equation}\label{II:eq:stages}
 K\Bigl(\frac\omega v\Bigr)^{k+1}\theta^m\le c^{zm}\,2^{-(k+1)}c^{-zm}<\frac12,\qquad
 \Bigl(\frac\omega v\Bigr)^j\theta^e<\frac12\quad(j\ge1,\ e\ge jT)
\end{equation}
(\lean{alg_total_volume}, \lean{alg_stage_lt_half}).

\begin{lemma}[\lean{exists_nalg_dimension}, \lean{exists_nalg_baseData}]\label{II:lem:data}
Let $k,z\ge1$, $t=z-1$, and let $\theta,A,B,c$ satisfy \eqref{II:eq:P} with $p=k+1$. For all
$a_0,g_0\ge0$, $\varepsilon'>0$, $\Gamma_0\in\R$ and $c_N>0$ there is $T\ge1$ such that
$c^{(k+1)T}\ge2$, the values \eqref{II:eq:values} are data for $k$ slopes, and
\[
 \frac{a_0}v+g_0\,Q\,B^T<\varepsilon',\qquad
 \Gamma_0<\frac{K\theta^m}{(m+1)\cdots(m+k+1)\,v^{k+1}A^m\,c_NQ^t} .
\]
\end{lemma}

\begin{proof}
For large $T$: the margins of the data follow from \eqref{II:eq:stages} and \cite[Lemma~8.3]{C},
as in \cite[Lemma~8.4]{C}; the torus margin holds because $(k+3)(m+k+2)/\sigma\le400(m+k+4)^3$
while $\omega/v=\frac12(c^z\theta)^{-T}$ grows exponentially ($c^z\theta<1$); and
$(1+\sigma)\theta<1$. Next $v=2(c^z\theta/B)^T\to\infty$ since $B<c^z\theta$, and
$QB^T\le(c^{k+1}B)^T\to0$ since $c^{k+1}B<1$. Finally $K\ge(c^m/2)^z$, $Q^t\le c^{mt}$ and
$v^{k+1}=2^{k+1}(c^z\theta/B)^m$, so the right side of the second inequality is at least
\[
 \frac{\varrho^{\,m}}{2^{z+k+1}c_N\,(m+1)\cdots(m+k+1)},\qquad\varrho=\frac B{Ac^t}>1,
\]
which tends to infinity, since $(m+1)\cdots(m+k+1)\le((2k+2)T)^{k+1}$.
\end{proof}

\begin{proposition}[Comparison; \lean{cofinal_alg_field_minors_inconsistent}]\label{II:prop:cmp}
Let $\nu>1$, let the setting of Section~\ref{II:sec:ar} be given, and let
$b_0,\dots,b_{K-1}\in(E^\times)^{k+1}$ be centers with a common denominator $d$ such that at
every embedding $\tau$ the points $\tau(b_c)$ have labelled logarithms \eqref{II:eq:centers}. If
\begin{equation}\label{II:eq:E}
 E_{\rm ar}(d)+E_{\rm an}<\mathrm{gap}\qquad\text{and}\qquad1+E_{\rm ar}(d)+E_{\rm an}<\Gamma,
\end{equation}
then it is impossible that for every $H_1$ there are an integer $H\ge H_1$ and a selection
with $\Delta_H\ne0$.
\end{proposition}

\begin{proof}
By Proposition~\ref{II:prop:upper} at every embedding and Corollary~\ref{II:cor:obs},
\[
 -(1-\overline b_H)-E_{\rm ar}(d)\le E_{\rm an}+\mathrm{rem}'_H
 +\max\bigl\{-\Gamma_H,\ -\nu\bigl(A(1-\eta)-\overline b_H\bigr)\bigr\}.
\]
If the maximum is the second term, then $\nu A(1-\eta)-1-(\nu-1)\overline b_H\le \mathcal E+\mathrm{rem}'_H$
with $\mathcal E=E_{\rm ar}(d)+E_{\rm an}$; the left side decreases in $\overline b_H\le\theta$, so
$\mathrm{gap}\le \mathcal E+\mathrm{rem}'_H$. Otherwise $\Gamma_H\le1+\mathcal E+\mathrm{rem}'_H$. Since
$\mathrm{rem}'_H\to0$ and $\Gamma_H\to\Gamma$, both contradict \eqref{II:eq:E} for large $H$
(\lean{cofinal_determinant_bounds_inconsistent}).
\end{proof}

\begin{proof}[Proof of Theorem~\ref{II:thm:nosupply}]
Suppose a $\Phi_k^{\mathcal H}$-supply of $\nu$-good vectors exists for every finite list
$\mathcal H$. Put $a=(a_0,\dots,a_{k-1},g)$ and $L=(L_0,\dots,L_{k-1},G)$; then $L_0\ne0$
(\lean{log_snoc_zero_ne}, using $k\ge1$) and $a$ has no power relation
(\lean{noPowerRelation_snoc_real}). The choices are made in the following order.
\begin{enumerate}
\item By Proposition~\ref{II:prop:family} the tuples $s_l(\alpha,\gamma)$ form a tuple family of
size $z=t+1$ for $a$, with a denominator $D$.
\item Choose $\theta,A,B,c$ by Lemma~\ref{II:lem:params} with $p=k+1$, $s=1$; choose
$\eta\in(0,1)$ with $\mathrm{gap}>0$ (\lean{exists_eta_preserving_exponent_gap}), and put
$\varepsilon=\mathrm{gap}/2$.
\item Let $N_b,c_R,c_N$ be the constants of Proposition~\ref{II:prop:supply},
$\Lambda_L=\sum_t|L_t|$, and
\[
 a_0=\Lambda_0F_0+\log2,\quad g_0=(k+1)z\log D+2\bigl(2c_R(\Lambda_L+N_b+1)+2\bigr),\quad
 \Gamma_0=\Bigl(\frac{1+\varepsilon}{c_0\,\eta^{(k+2)/(k+1)}}\Bigr)^{k+1}.
\]
\item Choose $T$ and the data by Lemma~\ref{II:lem:data} with $\varepsilon'=\varepsilon/2$; this
fixes $m,K,\omega,v,\sigma$ and $Q$. Let $\mathcal H_Q$ be the list of
Proposition~\ref{II:prop:supply}(T5$'$).
\item Put $R=2c_RQ(\Lambda_L+N_b+1)+2$ and choose $X$ with
$W/(X+1)<\varepsilon/2$, where $W=\log(2RM_L)+\nu+\log2+\log c_L$
(\lean{exists_weight_error_margin}); put $w_*=X+1$.
\item Let $T(\mathcal D)$ be the thresholds of Theorem~\ref{II:thm:U} for the data.
\item By Lemma~\ref{II:lem:sel}, applied to the $\Phi_k^{\mathcal H_Q}$-supply, choose $\nu$-good
vectors $r_0,\dots,r_{m-1}$ in general position and off $\mathcal H_Q$, with $q_i\ge2$, weights
above $T(\mathcal D)$ and $w_i>X+1$.
\item This gives the setting of Section~\ref{II:sec:ar} with $N=c_NQ^t$. The centers are the
$K=Q^z$ products $b_c$ of the tuple family. At every embedding they have labelled logarithms
\eqref{II:eq:centers} by Proposition~\ref{II:prop:supply}(T3,~T4$'$): $|s_cL_t|+\frac34\le
c_RQ\Lambda_L+\frac34\le R/2$ and $|Y^0_{\mu,t}|\le c_R(N_b+1)Q\le R$. Their common denominator is
$d=D^{zQ}$ by (T2).
\item The margins: $1/\omega=B^T$, so
\[
 E_{\rm ar}(d)+E_{\rm an}=\frac{a_0}v+\bigl((k+1)zQ\log D+2R\bigr)B^T+\frac W{w_*}
 \le\frac{a_0}v+g_0QB^T+\frac W{X+1}<\varepsilon<\mathrm{gap},
\]
and $\Gamma>c_0\eta^{(k+2)/(k+1)}\Gamma_0^{1/(k+1)}=1+\varepsilon>1+E_{\rm ar}(d)+E_{\rm an}$.
\item The $r_i$ are off $\mathcal H_Q$, so the centers have distinct characters at every
embedding. By Lemma~\ref{II:lem:minors} there are arbitrarily large $H$ with $\Delta_H\ne0$, which
contradicts Proposition~\ref{II:prop:cmp}.
\end{enumerate}
The list $\mathcal H_Q$ depends on $Q$, which is fixed in step 4, before the vectors are
selected in step 7.
\end{proof}

\section{The corollaries}\label{II:sec:cor}

\begin{proposition}[\lean{tupleIndep_card_le_finrank}]\label{II:prop:degree}
Let $E$ be a number field, $\alpha\in(E^\times)^n$ and $s_1,\dots,s_z$ an independent family
of automorphisms. Then $z\le[\Q(\alpha_1,\dots,\alpha_n):\Q]$.
\end{proposition}

\begin{proof}
If $s_i$ and $s_j$, $i\ne j$, agree on $F=\Q(\alpha)$, then $s_i(\alpha)s_j(\alpha)^{-1}=1$ is a
relation with $e=e_i-e_j\ne0$. So the restrictions $s_l|_F:F\to E$ are distinct
(\lean{tupleIndep_restrict_injective}). By Dedekind--Artin they are linearly independent in
$\operatorname{Hom}_\Q(F,E)$, which has dimension $[F:\Q]$ over $E$.
\end{proof}

\begin{proof}[Proof of Corollary~\ref{II:cor:main}]
(1) is Proposition~\ref{II:prop:degree} with Theorem~\ref{II:thm:main} (\lean{measureAlg_degree}).
In Lean the statement is also proved for real numbers given without a field
(\lean{measureAlg_real}): $E$ is the splitting field of $\prod_i\operatorname{minpoly}(a_i)$,
which is Galois, with an embedding into $\C$ that sends roots $\alpha_i$ to $a_i$; it maps
$\Q(\alpha)$ isomorphically onto $\Q(a)\subset\C$ (\lean{galois_model_tuple}). In the same
way Theorem~\ref{II:thm:main} is stated in Lean for real numbers given without a field
(\lean{measureAlg_real_rank}): a \emph{Galois model} of $a$ (\lean{GaloisModel}) is a Galois
number field with an embedding, a tuple $\alpha$ over $a$ and a largest independent family,
of size $z$; every tuple of positive real algebraic numbers has one
(\lean{exists_galoisModel}), and the bound holds with the exponent $n-1+z$.

(2) is (1) with $\Q(a)=\Q$ (\lean{measureAlg_rat}); (4) is (1) with
$\Q(a)\subset K$ (\lean{measureAlg_of_mem_field}, \lean{measureAlg_quadratic}); (3) is
Theorem~\ref{II:thm:main} for $n=2$: for two numbers the independence of conjugate tuples
is that of conjugate pairs (\lean{tupleIndep_pair_iff}), and $\Meas(2)$ gives
Corollary~\ref{I:cor:main}(5) with the hypotheses of Part~I
(\lean{twoLog_linearForm_of_measureAlg}).

(5) (\lean{measureAlg_quadUnits}). Take the Galois model above, so that $\alpha_i$ is rational or
satisfies $\alpha_i^2=t\alpha_i-e$ in $E$. For $g\in\Gal(E/\Q)$, $g(\alpha_i)$ is a root of the
same equation, whose roots are $\alpha_i$ and $t-\alpha_i=e\alpha_i^{-1}$; so
$g(\alpha_i)^2=\alpha_i^{\pm2}$ (\lean{conj_sq_of_quadUnit}). Let $s_1,\dots,s_z$ be independent
and $s_l(\alpha_i)^2=\alpha_i^{m_{l,i}}$. If $\sum_le_lm_{l,i}=0$ for all $i$, then
$\prod_ls_l(\alpha_i)^{2e_l}=\alpha_i^{\sum_le_lm_{l,i}}=1$, so $2e=0$. Hence the vectors
$(m_{l,i})_i\in\Z^n$ are linearly independent and $z\le n$ (\lean{tupleIndep_card_le_of_sq}).
\end{proof}

\begin{proof}[Proof of Example~\ref{II:ex:quad}]
Let $E$ be the compositum of the $s$ quadratic fields and $G=\Gal(E/\Q)$; let $\chi_i$ be the
character of $G$ by which it acts on the $i$-th field. Since $u_i$ is a unit of norm $\pm1$,
$g(u_i)=\pm u_i^{\chi_i(g)}$, and $g(p_j)=p_j$. The characters $\chi_i$ are distinct and
non-trivial. If $\prod_jp_j^{x_j}\prod_iu_i^{y_i}=1$, then
$\sum_jx_j\log p_j+\sum_iy_i\chi_i(g)\log u_i=0$ for all $g$; multiplying by $\chi_{i_0}(g)$
and summing over $G$ gives $y_{i_0}=0$, and then $x=0$. An element of $\Pi(a)$ is
$\prod_gg(a)^{x_g}$, with coordinates $\pm p_j^{y_0}$ and $\pm u_i^{y_i}$, where
$(y_0,y_1,\dots,y_s)=\sum_gx_g\bigl(1,\chi_1(g),\dots,\chi_s(g)\bigr)$. It is torsion only if
$y_i=0$ for all $i$ and, when $r\ge1$, $y_0=0$. The vectors $(1,\chi_1(g),\dots,\chi_s(g))$,
$g\in G$, span $\Q^{s+1}$, since the characters $1,\chi_1,\dots,\chi_s$ are linearly
independent; so $\rho(a)=s+1$ if $r\ge1$ and $\rho(a)=s$ if $r=0$. The field generated contains
$\Q(u_1,\dots,u_s)=E$, of degree $|G|=2^{s'}$ for some $s'\le s$, and $|G|\ge s+1$ because $G$ has the
$s+1$ distinct characters $1,\chi_1,\dots,\chi_s$.
\end{proof}

\part{Formal verification and remarks}

\section{Formal verification}\label{sec:lean}

\paragraph{Setting.}
The formalization is in Lean~4, version 4.34.1, with Mathlib at revision
\lean{d13f23b723b8a846827a245b89c10fc7d3f11612}, the revision used by the Lean files of \cite{C}; the
build report in the ancillary files records the toolchain, the package revisions and checksums
of the library.
It consists of 21 files \lean{QuadraticPi/TwoGrid*.lean} for Part~I (namespace
\lean{QuadraticPi.TwoGrid}) and 20 files \lean{QuadraticPi/NAlg*.lean} for Part~II (namespace
\lean{QuadraticPi.NAlg}), added to a project that contains the Lean supplement of \cite{L}
(which includes a selected part of the proof sources \cite{OAIL}), the files \lean{NLog*} and
\lean{CLog*} of \cite{C}; no existing file is modified. The new files accompany
this paper.

Part~I uses from the library of \cite{L} the two-torus interpolation theorem for
character-separated centers (Theorem~\ref{I:thm:interp},
\lean{TwoLogInterpolation.UniformStatement}), the truncation lemma, the weighted collision
bound with constant $\frac{\log2}8$, the lcm bound with $\Lambda=\log4+4$, the norm comparison
over a number field, the cofinal comparison of determinant bounds and the selection of
approximations above thresholds. The exponent $\mu$ is \lean{irrationalityExponent}, the
supremum of the $\nu>0$ with infinitely many $r\in\Q$, $\operatorname{den}r\ge2$,
$0<|x-r|<(\operatorname{den}r)^{-\nu}$; the library derives $2\le\mu(x)\le\kappa$ from the
finiteness of these sets for all $\nu>\kappa$.

Part~II uses the interpolation theorem of \cite{C} (Theorem~\ref{II:thm:U}), the homogeneous
selection, general position and conclusion lemmas, the truncation lemma, counting, the
collision and occupancy lemmas, the margins lemma and the arithmetic of entries of \cite{C},
the number field arithmetic and cofinal comparison of \cite{L}, and from Part~I only auxiliary definitions and lemmas: boxes, their indexing and torsion
coordinates (\lean{TwoGridPairsBasic}, and \lean{boxIdx}, \lean{boxCenter} of
\lean{TwoGridPairs}), the description of the embeddings of a Galois field
(\lean{TwoGridGalois}), two elementary inequalities of \lean{TwoGridParameters} and the
structure \lean{GridParams} of \lean{TwoGridData}. No theorem of Part~I on the irrationality
exponent is used in Part~II.

\paragraph{Main statements of Part~I.}
In Lean the Galois closure is an arbitrary Galois number field $E$ with an embedding
$\iota$, and $\iota(b_E)=b$, $\iota(a_E)=a$. The rank is the size $z$ of an independent
family $s$ of automorphisms (\lean{ConjIndep}) that is largest
(\lean{hmax}).
\begin{footnotesize}
\begin{verbatim}
theorem logRatio_irrationalityExponent_bounds [IsGalois ℚ E] (ι : E →+* ℂ)
    (bE aE : E) (b a : ℝ) (hb : 0 < b) (ha : 0 < a) (hιb : ι bE = b) (hιa : ι aE = a)
    (hmult : ∀ u v : ℤ, a ^ u * b ^ v = 1 → u = 0 ∧ v = 0)
    {z : ℕ} (s : Fin z → (E ≃ₐ[ℚ] E)) (hind : ConjIndep bE aE s)
    (hmax : ∀ (z' : ℕ) (s' : Fin z' → (E ≃ₐ[ℚ] E)), ConjIndep bE aE s' → z' ≤ z) :
    2 ≤ irrationalityExponent (Real.log a / Real.log b) ∧
      irrationalityExponent (Real.log a / Real.log b) ≤ (z : ℝ) + 2

theorem logRatio_degree_bounds_real (b a : ℝ) (hb : 0 < b) (ha : 0 < a)
    (hbi : IsIntegral ℚ b) (hai : IsIntegral ℚ a)
    (hmult : ∀ u v : ℤ, a ^ u * b ^ v = 1 → u = 0 ∧ v = 0) :
    2 ≤ irrationalityExponent (Real.log a / Real.log b) ∧
      irrationalityExponent (Real.log a / Real.log b) ≤
        (Module.finrank ℚ
          (IntermediateField.adjoin ℚ ({(b : ℂ), (a : ℂ)} : Set ℂ)) : ℝ) + 2

theorem logRatio_linearForm_bound_real (b a : ℝ) ... {ε : ℝ} (hε : 0 < ε) :
    ∃ c : ℝ, 0 < c ∧ ∀ p q : ℤ, (p ≠ 0 ∨ q ≠ 0) →
      c * (max 1 (max |(p : ℝ)| |(q : ℝ)|)) ^ (-((d : ℝ) + 1 + ε)) ≤
        |(p : ℝ) * Real.log a + q * Real.log b|
\end{verbatim}
\end{footnotesize}
In the third statement the hypotheses are those of the second, and $d$ stands for the
degree $[\Q(b,a):\Q]$, written out in Lean as in the second. The
first statement is Theorem~\ref{I:thm:main}, the second Corollary~\ref{I:cor:main}(1), the
third Corollary~\ref{I:cor:main}(5) in degree form; \lean{logRatio_linearForm_bound} is the
rank form.

\paragraph{Main statements of Part~II.}
The height \lean{height h} is $H(h)$, and \lean{NoPowerRelation} is the condition of the
introduction (both from \cite{C}).
\begin{footnotesize}
\begin{verbatim}
def TupleIndep {n : ℕ} (α : Fin n → E) {z : ℕ} (s : Fin z → (E ≃ₐ[ℚ] E)) : Prop :=
  ∀ e : Fin z → ℤ, (∀ i, ∏ l, (s l (α i)) ^ e l = 1) → e = 0

def MeasureAlg (n : ℕ) : Prop :=
  ∀ (E : Type) [Field E] [NumberField E] [IsGalois ℚ E] (ι : E →+* ℂ) (α : Fin n → E)
    (a : Fin n → ℝ), (∀ i, 0 < a i) → (∀ i, ι (α i) = a i) → NoPowerRelation a →
    ∀ (z : ℕ) (s : Fin z → (E ≃ₐ[ℚ] E)), TupleIndep α s →
      (∀ (z' : ℕ) (s' : Fin z' → (E ≃ₐ[ℚ] E)), TupleIndep α s' → z' ≤ z) →
      ∀ ε : ℝ, 0 < ε → ∃ C : ℝ, 0 < C ∧ ∀ h : Fin n → ℤ, h ≠ 0 →
        C * height h ^ (-((n : ℝ) - 1 + z) - ε) ≤ |∑ i, (h i : ℝ) * Real.log (a i)|

theorem measureAlg_all : ∀ n : ℕ, 1 ≤ n → MeasureAlg n

theorem tupleIndep_card_le_finrank {n : ℕ} (α : Fin n → E) (hα : ∀ i, α i ≠ 0) {z : ℕ}
    (s : Fin z → (E ≃ₐ[ℚ] E)) (hind : TupleIndep α s) :
    z ≤ Module.finrank ℚ (IntermediateField.adjoin ℚ (Set.range α))

theorem measureAlg_real {n : ℕ} (hn : 1 ≤ n) (a : Fin n → ℝ) (ha : ∀ i, 0 < a i)
    (hai : ∀ i, IsIntegral ℚ (a i)) (hind : NoPowerRelation a) (ε : ℝ) (hε : 0 < ε) :
    ∃ C : ℝ, 0 < C ∧ ∀ h : Fin n → ℤ, h ≠ 0 →
      C * height h ^ (-((n : ℝ) - 1 +
        Module.finrank ℚ (IntermediateField.adjoin ℚ (Set.range (fun i => (a i : ℂ))))) - ε) ≤
        |∑ i, (h i : ℝ) * Real.log (a i)|
\end{verbatim}
\end{footnotesize}
The first theorem is Theorem~\ref{II:thm:main}; in the second the field is written
\lean{tupleField} in the source, and the definitions are stated for a number field $E$; the
third is Corollary~\ref{II:cor:main}(1) for real numbers, and \lean{measureAlg_degree} is its
form for a given Galois field. The rank form for real numbers and the case $n=2$ are
\begin{footnotesize}
\begin{verbatim}
theorem measureAlg_real_rank {n : ℕ} (hn : 1 ≤ n) (a : Fin n → ℝ) (ha : ∀ i, 0 < a i)
    (hai : ∀ i, IsIntegral ℚ (a i)) (hind : NoPowerRelation a) (ε : ℝ) (hε : 0 < ε) :
    ∃ M : GaloisModel a, ∃ C : ℝ, 0 < C ∧ ∀ h : Fin n → ℤ, h ≠ 0 →
      C * height h ^ (-((n : ℝ) - 1 + M.z) - ε) ≤ |∑ i, (h i : ℝ) * Real.log (a i)|

theorem tupleIndep_pair_iff (bE aE : E) {z : ℕ} (s : Fin z → (E ≃ₐ[ℚ] E)) :
    TupleIndep (![aE, bE] : Fin 2 → E) s ↔ TwoGrid.ConjIndep bE aE s
\end{verbatim}
\end{footnotesize}
where \lean{GaloisModel a} is the structure of a Galois number field $E$, an embedding
$\iota$, a tuple $\alpha$ with $\iota(\alpha_i)=a_i$ and a largest independent family of
size \lean{z}. The bound holds in every Galois model, with the exponent $n-1+z$ of that model
(\lean{measureAlg_model}); \lean{measureAlg_real_rank} applies this to the model given by
\lean{exists_galoisModel}. The theorem \lean{twoLog_linearForm_of_measureAlg} has the
hypotheses and the conclusion of \lean{logRatio_linearForm_bound} of Part~I, for a field $E$ in
\lean{Type} (the universe of \lean{MeasureAlg}; Part~I allows any universe), and is proved
from \lean{measureAlg_all} for $n=2$, without the two-logarithm proof.

\paragraph{Correspondence, Part~I.}
\begin{center}\scriptsize
\begin{tabular}{@{}ll@{}}
\toprule
This paper & Lean declaration (file)\\
\midrule
Theorem~\ref{I:thm:main} & \lean{logRatio_irrationalityExponent_bounds},\\
 & \lean{logRatio_finite_good} (\lean{TwoGridGalois})\\
Corollary~\ref{I:cor:main}(1) & \lean{logRatio_degree_bounds}, \lean{conjIndep_card_le} (\lean{TwoGridDegree}),\\
 & \lean{logRatio_degree_bounds_real}, \lean{galois_model} (\lean{TwoGridDegreeReal})\\
Corollary~\ref{I:cor:main}(2) & \lean{quadratic_logRatio_bounds},\\
 & \lean{quadratic_normSigns_bounds} (\lean{TwoGridQuadratic})\\
 & \lean{log_two_div_log_goldenRatio_bounds} (\lean{TwoGridGolden})\\
Corollary~\ref{I:cor:main}(3) & \lean{ratOrNormOne_logRatio_bounds},\\
 & \lean{conj_sq_logRatio_bounds} (\lean{TwoGridNormOne})\\
 & \lean{log_one_add_sqrt_two_div_log_goldenRatio_bounds},\\
 & \lean{log_two_add_sqrt_three_div_log_one_add_sqrt_two_bounds},\\
 & \lean{quadUnits_eq_zpow} (\lean{TwoGridExamples})\\
Corollary~\ref{I:cor:main}(4) & \lean{log_cubicAlpha_div_log_two_bounds},\\
 & \lean{log_cubicAlpha_div_log_abs_bounds} (\lean{TwoGridCubic})\\
Lemma~\ref{I:lem:cubic} & \lean{cubicPoly_irreducible}, \lean{cubic_emb_cases},\\
 & \lean{log_indep_of_det} (\lean{TwoGridCubicField})\\
Corollary~\ref{I:cor:main}(5) & \lean{logRatio_linearForm_bound}, \lean{logRatio_linearForm_bound_real},\\
 & \lean{linearForm_lower_bound} (\lean{TwoGridLinearForms})\\
Lemma~\ref{I:lem:relation} & \lean{exists_relation_of_maximal},\\
 & \lean{exists_aut_of_embedding} (\lean{TwoGridGalois})\\
Proposition~\ref{I:prop:supply} & \lean{PairFamily}, \lean{PairFamily.centerSupply} (\lean{TwoGridPairs}),\\
 & \lean{exists_box_injective}, \lean{card_image_reduced_le} (\lean{TwoGridPairsBasic})\\
Definition~\ref{I:def:supply} & \lean{SupplyAt}, \lean{CenterSupply} (\lean{TwoGridTheorem})\\
Theorem~\ref{I:thm:B} & \lean{centerSupply_finite_good},\\
 & \lean{centerSupply_irrationalityExponent_bounds} (\lean{TwoGridTheorem})\\
Lemma~\ref{I:lem:params} & \lean{exists_twogrid_parameters} (\lean{TwoGridParameters})\\
choice of $m$ & \lean{exists_twogrid_dimension} (\lean{TwoGridDimension}),\\
 & \lean{exists_twogrid_baseData} (\lean{TwoGridBase})\\
matrix, non-zero minors & \lean{gridMatrix}, \lean{cofinal_field_grid_minors} (\lean{TwoGridMatrix})\\
Lemma~\ref{I:lem:arith} & \lean{gridSelectedMinor_arithmetic_lower_bound} (\lean{TwoGridArithmetic})\\
Lemma~\ref{I:lem:row} & \lean{period_coeff_eq_mainRowT} (\lean{TwoGridAnalyticSeries})\\
Proposition~\ref{I:prop:upper} & \lean{weighted_main_row_determinant_exp_boundT} (\lean{TwoGridAnalyticSeries}),\\
 & \lean{actualMinor_norm_upper_bound} (\lean{TwoGridAnalytic}), \lean{TwoGridData} (\lean{TwoGridData})\\
end of the proof & \lean{cofinal_grid_field_minors_inconsistent} (\lean{TwoGridArithmetic})\\
\bottomrule
\end{tabular}
\end{center}

\paragraph{Correspondence, Part~II.}
\begin{center}\scriptsize
\begin{tabular}{@{}ll@{}}
\toprule
This paper & Lean declaration (file)\\
\midrule
Theorem~\ref{II:thm:main} & \lean{measureAlg_all}, \lean{measureAlg_succ} (\lean{NAlgTheorem})\\
Corollary~\ref{II:cor:main}(1) & \lean{tupleIndep_card_le_finrank}, \lean{measureAlg_degree},\\
 & \lean{galois_model_tuple}, \lean{measureAlg_real} (\lean{NAlgDegree})\\
Theorem~\ref{II:thm:main} for real numbers & \lean{GaloisModel}, \lean{exists_galoisModel},\\
 & \lean{measureAlg_model}, \lean{measureAlg_real_rank} (\lean{NAlgBridge})\\
Corollary~\ref{II:cor:main}(3) & \lean{tupleIndep_pair_iff},\\
 & \lean{twoLog_linearForm_of_measureAlg} (\lean{NAlgBridge})\\
Corollary~\ref{II:cor:main}(2), (4) & \lean{measureAlg_rat}, \lean{measureAlg_of_mem_field},\\
 & \lean{measureAlg_quadratic} (\lean{NAlgExamples})\\
Corollary~\ref{II:cor:main}(5) & \lean{conj_sq_of_quadUnit}, \lean{tupleIndep_card_le_of_sq},\\
 & \lean{galois_model_apply}, \lean{measureAlg_quadUnits} (\lean{NAlgExamples})\\
$\Meas(n)$, independence, $\Meas(1)$ & \lean{MeasureAlg}, \lean{TupleIndep}, \lean{measureAlg_one} (\lean{NAlgTarget})\\
Lemma~\ref{II:lem:sel} & \lean{exists_general_position_above_thresholds_extra} (\lean{NAlgEndpoint})\\
Proposition~\ref{II:prop:bounded} & \lean{measureAlg_succ_of_bounded}, \lean{ordered_bound_shift} (\lean{NAlgConclusion})\\
Lemma~\ref{II:lem:mono} & \lean{tupleIndep_of_monomial} (\lean{NAlgTarget}),\\
 & \lean{tupleIndep_of_reduced} (\lean{NAlgHyperplane})\\
Lemma~\ref{II:lem:plane} & \lean{hyperplaneBounded_alg} (\lean{NAlgHyperplane})\\
Proposition~\ref{II:prop:nosupply} & \lean{measureAlg_succ_of_no_supply} (\lean{NAlgConclusion}),\\
 & \lean{supply_of_unbounded_of} (\lean{NAlgSelection})\\
Theorem~\ref{II:thm:nosupply} & \lean{no_supply_alg} (\lean{NAlgEndpoint})\\
Definition~\ref{II:def:family} & \lean{TupleFamily} (\lean{NAlgSupply})\\
Lemma~\ref{II:lem:relation}, Proposition~\ref{II:prop:family} & \lean{exists_relation_of_maximal_tuple}, \lean{exists_tupleFamily} (\lean{NAlgGalois})\\
Proposition~\ref{II:prop:supply} & \lean{center_integral}, \lean{exists_logs},\\
 & \lean{TupleFamily.exists_hyperplanes} (\lean{NAlgSupply})\\
Lemma~\ref{II:lem:T5} & \lean{exists_hyperplanes}, \lean{characterDistinct_box} (\lean{NAlgCharacters})\\
the setting & \lean{NAlgData} (\lean{NAlgData})\\
matrix, Lemma~\ref{II:lem:minors} & \lean{algMatrix}, \lean{cofinal_field_alg_minors} (\lean{NAlgMatrix})\\
Proposition~\ref{II:prop:lower}, Corollary~\ref{II:cor:obs} & \lean{algSelectedMinor_arithmetic_lower_bound},\\
 & \lean{algSelectedMinor_uniform_upper_bound_obstruction} (\lean{NAlgArithmetic})\\
\eqref{II:eq:centers}, Lemmas~\ref{II:lem:shift}, \ref{II:lem:det} & \lean{Centers}, \lean{mainEntry_eq_columnRow},\\
 & \lean{weighted_label_determinant_exp_bound} (\lean{NAlgAnalytic})\\
$\Gamma$, $E_{\rm an}$, low groups & \lean{collisionLimit}, \lean{analyticError}, \lean{low_label_capacity} (\lean{NAlgData})\\
Proposition~\ref{II:prop:upper} & \lean{actualMinor_norm_upper_bound} (\lean{NAlgAnalytic})\\
Lemma~\ref{II:lem:params} & \lean{exists_nalg_parameters} (\lean{NAlgParameters})\\
Lemma~\ref{II:lem:data} & \lean{exists_nalg_dimension}, \lean{exists_nalg_baseData} (\lean{NAlgDimension})\\
Proposition~\ref{II:prop:cmp} & \lean{cofinal_alg_field_minors_inconsistent} (\lean{NAlgComparison})\\
\bottomrule
\end{tabular}
\end{center}

\paragraph{Differences in presentation.}
In both parts the formal proof keeps the analytic estimate in the form of a remainder that
tends to $0$ along the degrees $H$, and compares the two sides through the cofinal comparison
lemma of \cite{L}. In Part~I, Corollary~\ref{I:cor:main}(1) for real $a,b$ is proved by
transferring any property of $[\Q(b_E,a_E):\Q]$ from the splitting field to
$\Q(a,b)\subset\C$ (\lean{galois_model}); Corollaries~\ref{I:cor:main}(3) and (4) are proved
with the pair families of Section~\ref{I:sec:cor} and Lemma~\ref{I:lem:cubic} directly, not
through the Galois closure; Corollary~\ref{I:cor:main}(2) with norms $\pm1$ is also the
theorem
\begin{center}\small\lean{quadraticNormSignsLogRatio_irrationalityExponent_bounds}\end{center}
of the file \lean{QuadraticNormOneLogRatioTheorem} of the library of \cite{C}.
In Part~II the matrix at the centers $b_c$ is the matrix of \cite{C} with $J=1$ and the base
$b_c$, row by row; so the truncation lemma, the integrality of the cleared entries and the
identification of the entries with coefficients of the chart are used as proved there.
Lemma~\ref{II:lem:det} is proved for an arbitrary grouping of the rows and an arbitrary shift
per group, with the shift entering the column function through $e^{h\cdot\mathrm{shift}}$.
Corollary~\ref{II:cor:main}(1) for real numbers is proved by transferring any property of
$[\Q(\alpha):\Q]$ from the splitting field to $\Q(a)\subset\C$.

\paragraph{Audit.}
The file \lean{NAlgAudit} prints the definitions and the main statements and runs
\texttt{\#print axioms} on 42 statements: for Part~II, Theorem~\ref{II:thm:main} with its
induction step and base, its forms for real numbers and in a Galois model, Corollary~\ref{II:cor:main}(1)--(5) and
the principal intermediate statements (Propositions~\ref{II:prop:bounded},
\ref{II:prop:nosupply}, \ref{II:prop:family}, \ref{II:prop:supply}, \ref{II:prop:lower},
\ref{II:prop:upper}, \ref{II:prop:cmp}, Lemmas~\ref{II:lem:plane}, \ref{II:lem:minors},
\ref{II:lem:params}, \ref{II:lem:data}, Theorem~\ref{II:thm:nosupply}); for Part~I,
Theorems~\ref{I:thm:main} and \ref{I:thm:B} and Corollary~\ref{I:cor:main}(1)--(5) (both
cases of (2)) with all the numerical examples. For each of them the axioms reported are \lean{propext},
\lean{Classical.choice} and \lean{Quot.sound}. The new files contain no \lean{sorry},
\lean{admit}, \lean{native_decide} or axiom declarations. Their SHA-256 checksums are:
\begin{center}\scriptsize
\begin{tabular}{@{}lll@{}}
\toprule
File & Lines & SHA-256\\
\midrule
\lean{TwoGridParameters} & 262 & \lean{b82660928fcb6bf799dfd5cf2777d079044b552fc1d06cff7d1f1d717ede8623}\\
\lean{TwoGridDimension} & 407 & \lean{ed40777568170cb1c71cbe7a8a29cb0b92b66d44d559924a156e5ed6fa05c79a}\\
\lean{TwoGridBase} & 187 & \lean{17ed19b821f645b689b05e17945ed6a3f0d500aa323bb4e449bb57ca0d0d625c}\\
\lean{TwoGridMatrix} & 253 & \lean{38f01dd4bbff2e034e58826ed91df1bf25bbbd848dda1e57478b33b12b686bb4}\\
\lean{TwoGridAnalyticSeries} & 280 & \lean{305b9178d0bb29fc2e23e7ea73d09ab798591c35ea29cf0b299a8b2ed55b90b0}\\
\lean{TwoGridData} & 349 & \lean{4e2aaed543c31f9fd88891cc54c6844937aa28f3c0f3959048abe8ad4d680647}\\
\lean{TwoGridAnalytic} & 351 & \lean{f30a8c1f28699916691116b45b18638315f8dafaed68ab38da6f6e8bf037f3c4}\\
\lean{TwoGridArithmetic} & 396 & \lean{c0c425113d889d9305b0ce47f60328d6a144763b5588731e7a60d433363b9e1f}\\
\lean{TwoGridTheorem} & 302 & \lean{0a8d703244e6e6cb66ae99238a605bfe724bcbcdf06e3cc41e02a385d0df523b}\\
\lean{TwoGridPairsBasic} & 398 & \lean{38ecd4fab1a3ba19fdda5922d31acc577bd35935885c6594111c64c437f8669f}\\
\lean{TwoGridPairs} & 295 & \lean{8c697d635c1cb11731e02f2295d1a906f11884454bc40ed439db1d6088d005b9}\\
\lean{TwoGridGalois} & 222 & \lean{b115c1b3b6dc2e9c1b4e0317d10fde86b83fdb2c0b0dac24c4457d87b46059db}\\
\lean{TwoGridQuadratic} & 274 & \lean{6b82d995e80c6e926a04d56f085b6a81b3d5cde3ca64ac176f07b8d7b32a8358}\\
\lean{TwoGridGolden} & 117 & \lean{1b8d44108ce03a3fb2d8c037ba93d37340420442b694bfe5c330a3d622675603}\\
\lean{TwoGridNormOne} & 181 & \lean{c21249a54512e20ebffe6c09fd1d1db34361fa8239589c5a0e061fbd12ae0709}\\
\lean{TwoGridExamples} & 273 & \lean{9da25d1389b6eb05ab11a98ae96b4d51108be2a4c96107820761d42937ca9dd8}\\
\lean{TwoGridCubicField} & 247 & \lean{cee2c345d236062e86392d9d69e87d5012721366ddc0717e69942c1a2fc9ccc6}\\
\lean{TwoGridCubic} & 303 & \lean{3f394197e695a6cc9676be7d44876bd626d72242b865f3159c1d82c244aa4ff5}\\
\lean{TwoGridDegree} & 102 & \lean{281e0cfa881337de89f0c3601f45a8f9562f166fa5f9f3abc26b67f7e9f98646}\\
\lean{TwoGridDegreeReal} & 95 & \lean{d56cf929bcb45874041012b7400533251de9cf5e993ea1e67721127d6f48981a}\\
\lean{TwoGridLinearForms} & 230 & \lean{4ecaf4b532dead4ea03262068edc28e124dc8565754274627b7f9173b6cc2d1b}\\
\bottomrule
\end{tabular}
\end{center}
\begin{center}\scriptsize
\begin{tabular}{@{}lll@{}}
\toprule
File & Lines & SHA-256\\
\midrule
\lean{NAlgTarget} & 140 & \lean{8a7614b9e312377e9dac1139846869ab2f06581cdf5fdc292e7609490df774f6}\\
\lean{NAlgParameters} & 209 & \lean{562b8f7efedeafcac09975d1d1f84abba235faeb4aea2b5dbe976a94863b01b0}\\
\lean{NAlgCharacters} & 229 & \lean{bd80a0076a3693607bd8d1a66bef5f6fcedea837d940407eeef58fdc34996b2f}\\
\lean{NAlgSupply} & 318 & \lean{f2024803a0d6aa3448827e5edeefe385faacdc33d3f262d350d18aa3208be275}\\
\lean{NAlgGalois} & 149 & \lean{9d55afbfe226e522f97745bfcccbdf84c5b333aa3bbb57831429971feb90f92d}\\
\lean{NAlgSelection} & 82 & \lean{2ff759d7d7e4805f5b4750b8e8f4ea07fd4e7e1a4c10a46efa8b26fe55bafb38}\\
\lean{NAlgHyperplane} & 287 & \lean{cab96eedb3c137ecdd7147702e14c05c0bf4b21548d2de8c1f3d4df3c7b4a55d}\\
\lean{NAlgConclusion} & 245 & \lean{8f1891fb80b98e1965459cf2a687a24deaec301807666ba8555be4367ef14103}\\
\lean{NAlgMatrix} & 255 & \lean{004125595c1dc822168a31b57e7ce334fbcfb543c2ae3dbb69e81c866aed3b16}\\
\lean{NAlgArithmetic} & 385 & \lean{f25a82ffdd8b0e662ada810245ba476b2d729a60f99548d4d9649764033bb478}\\
\lean{NAlgData} & 571 & \lean{39a848129fe4c7e98c6d5f5845c1485f45d41af6bf27d469f676c7b3de670f8d}\\
\lean{NAlgAnalytic} & 521 & \lean{e01a273adfa4526a30d1c5560aea7951608eae1e9865fdbe758efd3c22cfa0e3}\\
\lean{NAlgComparison} & 85 & \lean{38d7a6c81fecc47fc7155fa39db5f2bc47d3630b93248bbfcaa3dea861d16d52}\\
\lean{NAlgDimension} & 477 & \lean{7afe8d77db65b5532c4e85ec0363d3580f1928431987bf4fcbf021707e5d85bf}\\
\lean{NAlgEndpoint} & 394 & \lean{cc60e656ce45c3c91212a558f5367498c702b7503f141c3ae973a68442e74f62}\\
\lean{NAlgTheorem} & 35 & \lean{0b4d9ea20d8bfdc7a2d3c3db885c2b6e735df6ba317d8dd5b4324f9ba1303d99}\\
\lean{NAlgDegree} & 158 & \lean{179c7704ae0aefb47c1c77d54b9cd525aee427838f01884d2de3d4a640aad44c}\\
\lean{NAlgExamples} & 249 & \lean{fcee4de56dc3830ffa15559d857a02fec79cb31b1b137b6f6dd9961f63af5a3b}\\
\lean{NAlgBridge} & 173 & \lean{7922411bcc5ea9eda0a662e26b20841a52ca9aa45f2da092a49176aea4e1f8c8}\\
\lean{NAlgAudit} & 86 & \lean{2723573c3585535326c2c6a549496a115733e323d8ccf7438dccd127f6c5cea8}\\
\bottomrule
\end{tabular}
\end{center}

\paragraph{Not formalized.}
The identification of the rank $z$ of a Galois model with $\rk\Pi(a)$ for the Galois closure
in $\C$ (the rank of a group generated by a set is the largest size of an independent subset
of the set, and it does not depend on the Galois field containing the numbers); the passage
from Corollary~\ref{I:cor:main}(5), as derived from Part~II, back to the irrationality
exponent; the identities $\log\varphi=\artanh(1/\sqrt5)$ and $\log(1+\sqrt2)=\artanh(1/\sqrt2)$
in Corollary~\ref{I:cor:main}(3); the numerical check of Remark~\ref{I:rem:num};
Example~\ref{II:ex:quad} (the absence of a power relation and the computation of the rank); the necessity of $\nu>n+z$ in Remark~\ref{II:rem:where}; and the comparisons with
the literature in the introduction are not part of the formal verification.

\section{Remarks}\label{sec:remarks}

\begin{remark}[Part~I: where $\rk\Pi+2$ comes from]
For a family of rank $z$ the translation at $\tau$ is constant along the direction of
the pair sent to $(b,a)$; if this is the only direction killed at $\tau$, then $t=z-1$
needs every such pair in the span of the family, i.e.\ $z\ge\rk\Pi(a,b)$, and
$1+(t+2)/(z-t)=z+2$. When the pairs are matched ($\rk\Pi=1$) this is the geometric family
of \cite{L}, exponent $3$; for independent signs (Corollary~\ref{I:cor:main}(3) with two
different quadratic fields, or one rational number) $\rk\Pi=2$ and $4$ is the limit of
Theorem~\ref{I:thm:B}. This is a heuristic about the construction, not a theorem. Exponent
below $3$ would need removing the height condition $c^2B<1$.
\end{remark}

\begin{remark}[Part~II: where $n-1+\rho$ comes from]\label{II:rem:where}
With $K=Q^z$ centers and $N\asymp Q^{t}$ labels, the volume margin
$K(\omega/v)^{n}\theta^m<1$, the heights $QB^T\to0$ and the collision rate
$\Gamma\asymp\bigl(K\theta^m/(v^nA^mN)\bigr)^{1/n}\to\infty$ lead to the conditions
\eqref{II:eq:P} with $p=n$ and $s=z-t$. They are solvable if and only if $\nu s>z+p$, that
is, $\nu>1+(t+n)/(z-t)$. Sufficiency is Lemma~\ref{II:lem:params}. For the necessity,
\eqref{II:eq:P} gives $Ac^{t+p}<1$ and $c^s>A/\theta$, hence $A^{z+p}<\theta^{t+p}$, that is
$\frac{t+p}{z+p}<\frac{\log A}{\log\theta}$; since $-\log x/(1-x)$ decreases on $(0,1)$ and
$\theta<A$, $\frac{\log A}{\log\theta}<\frac{1-A}{1-\theta}$, and $\nu(A-\theta)>1-\theta$
says $\frac{1-A}{1-\theta}<1-\frac1\nu$; so $\frac{t+p}{z+p}<1-\frac1\nu$. A tuple family
has $t=z-1$, which gives $\nu>z+n$. For rational numbers $z=1$, $t=0$, $N=1$, and this is the
balance $\nu>n+1$ of \cite[Lemma~8.1]{C}; for $n=2$ it is the condition of
Lemma~\ref{I:lem:params}. A family with fewer labels at every embedding would
give a better exponent; for a tuple family of the largest size the relation of
Lemma~\ref{II:lem:relation} spends exactly one direction of the box.
\end{remark}

\begin{remark}[Constant terms]
For forms $h_0+h_1\log a_1+\dots+h_n\log a_n$ with algebraic $a_i$ the centers would need
additive coordinates, as in the second layer of \cite{C}; the single-logarithm case is
\cite{G}. This is not treated here.
\end{remark}

\begin{remark}[Ineffectivity]
The constants are not effective. In Part~II, for each $n\ge2$ the proof shows that there
cannot be $m$ good vectors in general position with denominators growing as fast as the
thresholds require; it gives no bound for the denominators themselves, and the induction uses
the constants of $\Meas(n-1)$ in Lemma~\ref{II:lem:plane}.
\end{remark}

\end{document}